\documentclass[letterpaper]{siamltex}
\usepackage[ansinew,latin1]{inputenc}
\usepackage[francais,english]{babel}
\usepackage{amsfonts,amssymb}        %
\usepackage{amsmath}
\usepackage{latexsym}
\usepackage{verbatim}    %
\usepackage{float}
\usepackage{graphicx}
\usepackage{wrapfig}
\usepackage{sidecap}
\usepackage{subfigure}
\def\ZETAM{L}
\def\KAPPAM{K}
\def\ID{{\mathrm{Id}\,}}
\def\INTP{\mathrm{int}}
\def\EXTP{\mathrm{ext}}

\def\TAG#1{{\qquad ({#1})}}
\def\TAGG#1{{\qquad ({#1})}}

\def\BIGO{\mathcal{O}}
\def\VIZ#1{(\ref{#1})}

\def\WEAKLY{{\rightharpoonup}}

\def\HESS{\mathrm{Hess}\,}
\def\SEP{\,|\,}

\def\DET{\mathrm{det}\,}
\def\TR{\mathrm{Tr}\,}

\def\DIAG{\mathrm{diag}\,}

\def\diss{\M{d}}
\def\nc{n}
\def\SMAN{{\mathcal{M}_{0}}}
\def\SMANZ#1{\mathcal{M}_{#1}}
\def\MASST{M}
\def\MASSTZ{M_z}

\def\nueps{\nu_{\epsilon}}
\def\pen#1{#1_{\nu}}
\def\nusc{\bar{\nu}}
\def\peneps#1{#1_{\nueps}}

\def\DIAG{\mathrm{diag}\,}

\def\Dp{\nabla_{\! p}}
\def\Dq{ \nabla_{\! q} }
\def\Dz{\nabla_{\! z}}
\def\Done{ \nabla_{\! 1} }
\def\Dtwo{\nabla_{\! 2}}

\def\penN#1{#1_{\nu_{N}}}

\def\dt{\delta t}
\def\dg{\nabla^{d}}
\def\dl{\Delta_{d}}

\def\EXPECT#1{\E\left[{#1}\right]}

\newcommand{\poisson}[1]{\left\{#1\right\}}

\def\E{\mathbb{E}}

\def\R{\mathbb{R}}

\def\EXP#1{e^{#1}}
\def\ds{\displaystyle}
\def\COMMA{\,,}
\def\PERIOD{\,.}
\def\M#1{{\rm #1 }}

\newcommand{\ph}{\varphi}

\newcommand{\norm}[1]{\left\Vert#1 \right\Vert}

\newcommand{\normop}[1]{|\!|\!|  #1 |\!|\!| }

\newcommand{\abs}[1]{\left\vert#1\right\vert}
\newcommand{\set}[1]{\left\{#1\right\}}

\newcommand{\limop}[1]{\mathop{#1}\limits}

\newenvironment{system} { \left \{
  \begin{alignedat}{2} }
{  \end{alignedat}
  \right.}

\newcommand{\syst}[1]{ \begin{system} #1 \end{system}}

\newcommand{\pare}[1]{\left(#1\right)}

\newenvironment{li} { \left .
  \begin{array}{lcl} }
{  \end{array}
  \right .}
\newcommand{\arr}[1]{ \begin{li} #1 \end{li}}

\newenvironment{tb}{\begin{tabular}{lcl}}{\end{tabular}}

  \newcommand{\bmat}{\begin{pmatrix}}
\newcommand{\emat}{\end{pmatrix}}

\newtheorem{The}{Theorem}[section]
\newtheorem{Lem}[The]{Lemma}
\newtheorem{Pro}[The]{Proposition}

\newtheorem{Def}[The]{Definition}

\newtheorem{scheme}{Scheme}[section]

\newtheorem{Rem}[The]{Remark}
\newtheorem{Exa}[The]{Example}

\newcommand{\bi}{\begin{itemize}}
\newcommand{\ei}{\end{itemize}}

\newcommand{\be}{\begin{equation}}
\newcommand{\ee}{\end{equation}}
\newcommand{\bea}{\begin{eqnarray}}
\newcommand{\eea}{\end{eqnarray}}

\newcommand{\bpro}{\begin{Pro}}
\newcommand{\epro}{\end{Pro}}
\newcommand{\blem}{\begin{Lem}}
\newcommand{\elem}{\end{Lem}}
\newcommand{\bthe}{\begin{The}}
\newcommand{\ethe}{\end{The}}

\newcommand{\bdfn}{\begin{Def}}
\newcommand{\edfn}{\end{Def}}
\newcommand{\brem}{\begin{Rem}}
\newcommand{\erem}{\end{Rem}}
\newcommand{\bexa}{\begin{Exa}}
\newcommand{\eexa}{\end{Exa}}

\title{Exact and non-stiff sampling of highly oscillatory systems:\\
       an implicit mass-matrix penalization approach}

\author{Petr Plech\'a\v{c}\thanks{%
               Department of Mathematics, University of Tennessee and Oak Ridge National Laboratory,
               Knoxville, TN 37996-1300, USA,
              ({\tt plechac@math.utk.edu}).
	}
	\and
	Mathias Rousset\thanks{%
               INRIA Lille Nord-Europe, Villeneuve d'Ascq, France, ({\tt mathias.rousset@inria.fr})}
}

\begin{document}
\maketitle

\begin{abstract}
We propose and analyze an implicit mass-matrix penalization (IMMP) technique which enables efficient and exact sampling
of the (Boltzmann/Gibbs) canonical distribution associated to Hamiltonian systems with fast degrees of freedom (fDOFs).
The penalty parameters enable arbitrary tuning of the timescale for the selected fDOFs, and the method is interpreted as
an interpolation between the exact Hamiltonian dynamics and the dynamics with infinitely slow fDOFs (equivalent to
geometrically corrected rigid constraints). This property translates in the associated numerical methods into a tunable trade-off between stability and dynamical modification.
The penalization is based on an extended Hamiltonian with artificial constraints associated with each fDOF.  By construction, the resulting dynamics is
statistically exact with respect to the canonical distribution in position variables.

The algorithms can be easily implemented with standard geometric integrators with algebraic constraints given by the expected fDOFs,
and has no additional complexity in terms of enforcing the constraint and force evaluations.
The method is demonstrated on a high dimensional system with non-convex interactions. Prescribing the macroscopic dynamical timescale, it is shown that the IMMP method increases the time-step stability region with a gain that grows
linearly with the size of the system. The latter property, as well as consistency of the macroscopic dynamics of the IMMP method is proved rigorously for linear interactions.
Finally, when a large stiffness parameter is introduced, the IMMP method can be tuned to be asymptotically stable,
converging towards the heuristically expected Markovian effective dynamics on the slow manifold.
\end{abstract}



\pagestyle{myheadings}
\thispagestyle{plain}
\markboth{P. Plech\'a\v{c}, M. Rousset}%
         {Sampling highly oscillatory systems.}

\section{Introduction}\label{s:intro}

This paper deals with stochastic numerical integration and sampling of Hamiltonian
systems with multiple timescales. The main motivation is to develop numerical methods which sample
{\it exactly} canonical distributions, while integrating the fastest degrees of freedom only statistically.
Furthermore, one also seeks good approximation of dynamical behavior, at least at large temporal scales.

Hamiltonian systems with multiple timescales typically appear in molecular dynamics (MD)
simulations, which have become, with the aid of increasing computational power, a standard tool
in many fields of physics, chemistry and biology.
However, extending the simulations to physically relevant time-scales remains
a major challenge for various large molecular systems.
Due to the complexity of implicit methods, the time scales reachable by standard numerical methods
are usually limited by the rapid oscillations of some particular degrees of freedom. Since the sampling dynamics has
to be integrated for long times, the time-step restriction associated with fast oscillations/short
time scales in molecular systems contributes to the high computational cost of such methods.
Yet  the physical necessity of resolving the fast degrees of freedom in simulations is often ambiguous, and efficient
treatment of the fast time scales has motivated new interest in developing numerical schemes
for the integration of such stiff systems.

The problem of integrating stiff forces is relevant both for the direct numerical simulation of the Hamiltonian dynamics,
and for the less restrictive problem of designing a process that samples the canonical ensemble.
Sampling from the canonical distribution can be achieved by Markov chain Monte Carlo (MCMC) algorithms based on an a priori
knowledge of possible moves, with a Metropolis-Hastings acceptance/rejection corrector (a historical reference is \cite{Met53}). For complex molecular systems however,
such global moves remain unknown in general, and sampling methods consists generically in using either an Hamiltonian dynamics with a thermostat (e.g., a Langevin process),
or a drifted random walk (Brownian dynamics), corresponding to over-damped Hamiltonian dynamics (see \cite{CanLeg05}
for a review and references on classical sampling methods). Brownian dynamics of systems with multiple time scales suffers
from similar stability restrictions (see \cite{Fix86,GunBer81} for some practical issues on Brownian dynamics simulation for
Molecular Dynamics).

Broadly speaking, one may start by recognizing two approaches to the numerical treatment of stiff systems:
\begin{description}
\item[{\rm (i)}]   Semi-implicit, multi-step integrators and their variants (e.g., the textbooks
     \cite{LeiRei05, HaiLub02}, or the review paper \cite{JahLub06} and references therein, \cite{Fix86} for Brownian dynamics),
     which attempt to resolve microscopic highly oscillatory dynamical behavior.
\item[{\rm (ii)}] Methods with direct constraints, where the highly oscillatory
     degrees of freedom are constrained to their equilibrium value
     (e.g., \cite{LeiRei05, LeiSke94, GunBer77, Fix78} and references therein).
\end{description}
In spite of their differences the common key feature of all these methods is to balance a trade-off between stability restrictions and implicit time-stepping form, or in other words,
between the computational effort associated with small time steps, and the computational cost of solving implicit equations implied by the stiffness.

The method developed in this paper aims at giving an appropriate interpolation between exact dynamics and
constrained dynamics considered in the second family of methods mentioned above. Although constrained dynamics remove,
in principle, the stiffness of the associated numerical scheme, it introduces new difficulties and numerical problems.
As an approximation to the original dynamics it modifies important features of the system; most importantly,
the original statistical distribution.
The principal goal of the proposed method is to replace direct constraints by implicit mass-matrix
penalization (IMMP), detailed in Section~\ref{s:IMMP},  which integrates fDOFs, but with a tunable mass penalty.
This approach guarantees that the canonical distribution is computed {\it exactly},
and that a freely tunable trade-off between dynamical modification and stability can be obtained.

The idea of adjusting mass tensors in order to slow down fast degrees of freedom goes back to \cite{Ben75}.
In this paper, the author proposes to modify the mass tensor with respect to the Hessian of the potential energy function
in order to confine the frequency spectrum to low frequencies only. Two natural drawbacks of this procedure arise
from the costly computation of the second-order derivatives of the potential, and from the bias introduced when the adjusted mass-tensor is adapted during the dynamics. Such an approach seems inevitable when the fDOFS are unknown, but in many cases, the fast degrees of freedom are explicitly given by the structure of the system (e.g., co-valent and angle bonds in molecular chains).
To our knowledge, mass tensor modification have been used in practical MD simulations by increasing the mass of some well-chosen (e.g., light) atoms \cite{FeeHesBer90,MaoFrie90}. The aim of this paper is to propose a more systematic mass-tensor modification strategy.

The proposed method relies on the assumption that the system Hamiltonian is separable with quadratic kinetic energy
\be\label{e:H}
H(p,q) = \frac{1}{2} p^{T} \MASST^{-1} p  + V(q) \COMMA
\ee
and that the ``fast'' degrees of freedom $(\xi_1,..,\xi_n)$ are explicitly defined smooth functions of the system position
\begin{equation}\label{e:xi}
q=(q_1,\dots,q_d) \mapsto \pare{\xi_{1}(q),...,\xi_{\nc}(q)}\PERIOD
\end{equation}
We emphasize that the knowledge of  ``fast forces'' is not required, and the variables $\xi$ can be chosen arbitrarily.
If the fDOFs are not identified the method retains its approximation properties while not performing efficiently.
The fDOFs are penalized with a mass-tensor modification given by
\be\label{e:PMM}
	\pen{\MASST}(q) = \MASST + \nu^{2}\nabla_{q}\xi \,\MASSTZ\,  \nabla_{q}^T\xi \COMMA
\ee
where $\nu$ denotes the penalty intensity, and $\MASSTZ$ a ``virtual'' mass matrix associated with the fDOFs. The modification does not impact motions orthogonal to the fDOFs.
The position dependence of the mass-penalization introduces a geometric bias. This bias is corrected by introducing an
effective potential
\be\label{e:Vfixnu}
V_{\M{fix},\nu}(q) = \frac{1}{2\beta} \ln \pare{ \DET (\pen{\MASST}(q)) } \COMMA
\ee
which will turn out to be a $\tfrac{1}{\nu}$-perturbation of the usual Fixman corrector (see \cite{Fix78}) associated with the sub-manifolds defined by constraining the fDOFs $\xi$.
The key point is then to use an implicit representation
of the mass penalty with the aid of the extended Hamiltonian
\begin{equation}\label{e:Himmp}
\syst{
& H_{\M{IMMP}}(p,p_{z},q,z) = \frac{1}{2} p^{T} \MASST^{-1} p  + \frac{1}{2} p_{z}^{T} \MASSTZ^{-1} p_{z}
             + V(q) + V_{\M{fix},\nu}(q)\COMMA & \\
& \xi(q) = \frac{z}{\nu} \PERIOD & \TAG{\pen{C}}
}\end{equation}
The auxiliary degrees of freedom $z$ are
endowed with the ``virtual'' mass-matrix $\MASSTZ$. The constraints $(\pen{C})$ are applied in order to identify
the auxiliary variables and the fDOFs $\xi$ with a coupling intensity tuned by $\nu$. The typical time scale of the fDOFs is thus enforced by the penalty $\nu$.  The system is coupled to a thermostat through a Langevin equation \eqref{e:langevin}, which yields a stochastically perturbed dynamics that samples the equilibrium canonical distribution. We then obtain the following desirable properties:
 \begin{enumerate}
 \item The associated canonical equilibrium distribution in position is independent of the penalty $\nu$.
 \item The limit of vanishing penalization ($\nu = 0$) is the original full dynamics,  enabling the construction of dynamically consistent numerical schemes.
 \item  The limit of infinite penalization is a standard effective constrained dynamics on the ''slow'' manifold associated with stiff constraints on $\xi$.
 \item Numerical integrators can be obtained through a simple modification of standard integrators for effective dynamics with constraints, with equivalent computational complexity.
 \end{enumerate}

 The dynamics associated with the IMMP Hamiltonian \eqref{e:Himmp} is detailed in \eqref{e:IMMP}, see Section~\ref{s:IMMP}.
Its numerical discretization by a leapfrog/Verlet scheme with constraints (usually called ``RATTLE'')
is given by \eqref{e:immpscheme}. The flow can then be corrected by a Metropolis step to obtain exact canonical sampling,
which is reminiscent of so-called Hybrid Monte-Carlo methods (see references in Section~\ref{s:numerinteg}).
When this correction is introduced, the gradient of the Fixman potential \eqref{e:Vfixnu} need not be computed.
The numerical aspects of the method are presented in Section~\ref{s:numerinteg}.

By including a penalty, the proposed method modifies the original Hamiltonian.
However, the mass penalty can be also thought of as depending on the time step
$\nu := \nu(\dt)$ leading to consistent schemes. In Section~\ref{s:Nanal} and
in Section~\ref{s:highoscill}, the penalty (and the fastest timescales) will even grow to infinity
$\nu \to + \infty$ as the dimension of the system $N\to\infty$  or a stiffness parameter $\epsilon\to 0$.
In both cases, we encounter generic situations where the IMMP method leads to dynamically consistent
limits, a partial differential equation, and an effective dynamics on a ``slow'' manifold, respectively.

High-dimensional systems usually contain a large variety of timescales, and are therefore challenging test cases.
A linear chain of particles with repulsive {\it non-convex} interactions is numerically studied in Section~\ref{s:numerexperiment}.
It is shown that a large mass-penalty of order $\nu :=\nusc N $, where $N$ is the system size, induces a gain in time stepping
of order $N$, while numerical evidence are given that the macroscopic dynamics, given by a formal stochastic partial differential equation,
remains of order $1$ (in particular, the convergence to equilibrium). Rigorous proofs with explicit scalings of these behaviors is provided for the linear case (i.e. harmonic interactions) in Section~\ref{s:Nanal}, and consistence of the IMMP macroscopic dynamics when the scaled penalty vanishes $\nusc \to 0$ is even demonstrated.

In Section~\ref{s:highoscill},  we introduce the stiffness  parameter $\epsilon$ and show that the penalty intensity can be scaled as $\nu:= \nusc/\epsilon$ in order to obtain asymptotically stable dynamics in the limit $\epsilon \to 0$. We prove that the dynamics converge towards the expected Markovian effective dynamics on the slow manifold. We also present analysis of stability properties of the proposed scheme.

\medskip\noindent{\bf Acknowledgments:}
The research of M.R. was partially supported by the EPSRC grant GR/S70883/01 while he was visiting Mathematics Institute,
University of Warwick. The research of P.P. was partially supported by the Office of Advanced Scientific Computing Research,
U.S. Department of Energy; the work was partly done at the ORNL, which is managed by UT-Battelle, LLC under Contract No. DE-AC05-00OR22725.

\section{Formulation}\label{s:notations}
We consider a Hamiltonian system in the phase-space $\R^{d}\times\R^{d}$
with the Hamiltonian $H$ in the form
\be\label{e:HH}
H(p,q) = \frac{1}{2} p^{T} \MASST^{-1} p  + V(q) \COMMA
\ee
We use generic matrix notation, for instance the Euclidean scalar product of two vectors $p_1,p_2\in\R^N$
denoted by $p_1^T p_2$, and the gradients of mappings from $\R^d$ to $\R^n$ with respect to standard bases are represented by matrices
\[
 (\nabla^T_q\xi)_{ij} = (\nabla_q\xi)_{ji} = \frac{\partial \xi_i}{\partial q_j}\COMMA\; i=1,\dots,n\,,\;
                                                                                         j=1,\dots,d \PERIOD
\]

When the system is thermostatted, i.e., kept at the constant temperature, the long time distribution of the system in the
phase-space is given by the canonical equilibrium measure at the inverse temperature $\beta$
(also called the NVT distribution) given by
\be \label{e:boltzmann}
\mu(dp\,dq) = \frac{1}{Z}\EXP{ - \beta H(p,q) } dp\, dq \COMMA\;\;\; Z = \int_{\R^d\times\R^d} \EXP{ - \beta H(p,q) } dp\,dq\COMMA
\ee
with the normalization constant $Z<\infty$. The standard dynamics used to model thermostatted systems
are given by Langevin processes.
\begin{Def}[Langevin process]\label{d:langevin} A Langevin process at the inverse temperature $\beta$ with the Hamiltonian
  $H(q,p)$, $(p,q)\in\R^d\times \R^d$, the $d\times d$ dissipation matrix $\gamma$, and the fluctuation matrix $\sigma$
  is given by the stochastic differential equations
  \be\label{e:langevin}
  \syst{
    &\dot{q} = \Dp H  &\\
    &\dot{p} = -\Dq H -\gamma \dot{q} + \sigma \dot{W} \COMMA &
    }
   \ee
  where $\dot{W}$\footnote{Throughout the paper, stochastic integrands have finite variation,
  so that stochastic integration (e.g. Ito or Stratonovitch) need not be specified.} is a standard white noise (Wiener process), and
  $\sigma\in\R^d\times\R^d$
  satisfies the fluctuation-dissipation identity
  \[
  \sigma \sigma^{T} =\frac{2}{\beta} \gamma \PERIOD
  \]
For any $\gamma$, the process is reversible with respect to the stationary canonical distribution \VIZ{e:boltzmann}.
Furthermore, if $\gamma$ is strictly positive definite, the process  is ergodic.
\end{Def}
Throughout this paper, the usual global
Lipschitz conditions  (see \cite{Osk92}) on $H$ and $\xi$ are assumed, ensuring well-posedness of the considered stochastic differential equations. The analysis presented in the paper can be generalized to a position dependent
dissipation matrix $\gamma = \gamma(q)$.

The mapping $\xi:\R^d\to\R^n$, defines
$\nc \leq d$ degrees of freedom, given by smooth  functions
taking values in a neighborhood of $0$. We assume that the mapping $\xi$ is regular
(i.e., with a non-degenerate Jacobian) in an open $\delta$-neighborhood
$\mathcal{O}_{\delta}=\set{q\SEP \norm{\xi(q)} < \delta} $ of $\xi^{-1}(0)$, hence defining a smooth sub-manifolds of
$\mathbb{R}^{d}$ denoted $\SMANZ{z}=\xi^{-1}(z)$ for $z$ in a neighborhood of $0$.
The dependence of the potential with respect to the degrees of freedom $\xi$ is expected to be ``stiff'' in the second variable. In Section \ref{s:highoscill} we will introduce the stiffness parameter
$\epsilon$. We will assume in this section that such parameter dependence can be explicitly identified, and that the potential energy
$V$ can be written in the form
\begin{equation}\label{Voscill}
V(q)= U(q,\frac{\xi(q)}{\epsilon}) \COMMA
\end{equation}
where the function $U:\R^d\times\R^n\to\R$ satisfies the coercivity condition
$\lim_{z \to +\infty} U(q,z) = +\infty$. The fast degrees of freedom $\xi$
of states at a given energy then remain in a closed neighborhood of the origin as the
stiffness parameter $\epsilon \to 0$. In this limit the system is confined to the sub-manifold $\SMAN$,
which is usually called the ``slow manifold''.

\section{The implicit mass-matrix penalization method}\label{s:IMMP}
In this section we focus on properties of the IMMP method. The multiscale
structure of the potential $V$ need not be known in order to apply the method.
Thus, in this section, we consider the potential $V$ in the form where we do not impose the
structural assumption \VIZ{Voscill} on the potential function $V:\R^d\to\R$.

The new, penalized mass-matrix of the system is the position dependent tensor defined in \eqref{e:PMM}. The associated modified impulses are denoted
\be\label{e:ppen}
	\pen{p} = \pen{\MASST}(q) \MASST^{-1} p \PERIOD
\ee
When $\nu$ becomes large, the velocities are bound to remain tangent to the manifolds $\set{q|\xi(q) = z }$,
and orthogonal motions are arbitrarily
slown down. Conversely, when $\nu=0$, one recovers the original highly oscillatory system.
Since the modification in $\pen{\MASST}$ depends on the position $q$,  new geometry is introduced and an additional correction \eqref{e:Vfixnu} in the potential energy is required in order to preserve original statistics in the position variable. This correction is in fact close to the standard Fixman corrector  for $\nu$ large (see \eqref{e:Vfix}).
Defining $G(q)$ as the $\nc \times \nc$ Gram matrix associated with the fast degrees of freedom
\begin{equation}\label{e:Gram}
G(q) =   \nabla^{T}_{q} \xi \,\MASST^{-1}  \Dq \xi \COMMA
\end{equation}
one has the following property of the correcting potential.
\begin{Pro} Up to an additive constant, we have
\be\label{e:penV}
V_{\M{fix},\nu}(q) =  \frac{1}{2\beta} \ln \DET \pare{ G(q) +   \frac{1}{\nu^{2} } \MASSTZ^{-1} } \COMMA
\ee
and thus (up to additive constants)
\[
\lim_{\nu \to +\infty} V_{\M{fix},\nu} = V_{\M{fix}} = \frac{1}{2\beta} \ln \DET \pare{ G(q) } \COMMA\;\;\;\;\mbox{and }\;\;\;\;\;
\lim_{\nu \to 0} V_{\M{fix},\nu} = 0 \PERIOD
\]
\end{Pro}
\begin{proof}
Using the  identity for a non-diagonal matrix $J$ of dimension $n_{1} \times n_{2}$:
\[
\DET( \ID_{n_1}+ JJ^{T}   ) = \DET( \ID_{n_2}+ J^{T}J   ) \COMMA
\]
one observes
\[
\DET(\pen{\MASST}) = \DET(\MASST)\, \DET(\nu ^{2}\MASSTZ)\, \DET(G+\frac{1}{\nu^2}\MASSTZ^{-1})
\]
from which the expression for the corrected Fixman potential follows.
\end{proof}

The associated modified Hamiltonian is then given by
\be\label{e:Hnu}
	\pen{H}(\pen{p},q) = \frac{1}{2} \pen{p}^{T} {\pen{\MASST}}^{-1} \pen{p}    + V(q) + V_{\M{fix},\nu}(q)  \COMMA
\ee
and $H_{0}=H$ is the original Hamiltonian \eqref{e:H}.

By construction, statistics of positions $q$ of the mass penalized Hamiltonian are independent of the penalization,
leading to the \emph{exact canonical statistics} in  position variables.
\bpro[Exact statistics]
The canonical distribution associated with the mass-penalized Hamiltonian \eqref{e:Hnu} is given by
\be\label{e:penboltz}
\pen{\mu}(d\pen{p} \, dq)= \frac{1}{\pen{Z}} \EXP{ - \beta \pen{H}(\pen{p},q) } d\pen{p}\,dq \PERIOD
\ee
Its marginal probability distribution in $q$ is
\[
 \frac{1}{\pen{Z}} \int  \EXP{ - \beta \pen{H}(\pen{p},q) } d\pen{p} = \frac{\EXP{ - \beta V(q) } dq}{\int \EXP{ - \beta V(q) }\,dq}
\]
which is the original canonical distribution \eqref{e:boltzmann} in the position variables,
and is independent of the mass penalization parameter $\nu$.
\epro
\begin{proof}
The normalization of Gaussian integrals in the $\pen{p}$ variables yields
\[
\int \EXP{ - \beta \frac{1}{2} \pen{p}^{T} {\pen{\MASST}}^{-1} \pen{p} } d \pen{p} =
 \left(\frac{2 \pi}{\beta}\right)^{d/2}\!\sqrt{\DET(\pen{\MASST})} \COMMA
\]
which is cancelled out  by the Fixman corrector $V_{\M{fix},\nu}$ and the result follows.
\end{proof}

Sampling such a system can be done using the standard Langevin stochastic perturbation
as detailed in Definition~\ref{d:langevin}.
However, the direct discretization of the equation of motion given by $\pen{H}$ (e.g., by an explicit scheme)
is bound to be unstable (from non-linear instabilities) when the fast degrees of freedom are not affine functions.
In order to construct stable schemes one may rather use an implicit formulation of the Hamiltonian \eqref{e:Hnu},
in conjunction with a solver which enforces the constraints.
To obtain such a formulation we extend the state space with $\nc$ new variables $(z_{1},..,z_{\nc})$,
and associated moments $(p_{z_{1}},..,p_{z_{\nc}})$. The auxiliary mass-matrix for the
new degrees of freedom is then given by $\MASSTZ$.
The new extended Hamiltonian of the system $H_{\M{IMMP}}$, defined by \eqref{e:Himmp}, is now defined in $\mathbb{R}^{d+\nc}\times\R^{d+\nc}$, where $\nc$ position constraints denoted by $(\pen{C})$ are included.
This construction implies $\nc$ hidden constraints on momenta. The equivalence of the two Hamiltonians
\eqref{e:Hnu} and  \eqref{e:Himmp} formulations is stated as a simple separate lemma.
\begin{Lem}
The equations of motion associated with the penalized mass-matrix Hamiltonian \eqref{e:Hnu} or the extended Hamiltonian
with constraints \eqref{e:Himmp} are identical.
\end{Lem}
\begin{proof}
The Lagrangian associated with $H_{\M{IMMP}}$ is given by
 \[
 L_{\M{IMMP}}(\dot{q},\dot{z},q,z) = \frac{1}{2} \dot{q}^{T} \MASST \dot{q}  +\frac{1}{2} \dot{z}^{T} \MASSTZ \dot{z}
- V(q) -  V_{\M{fix},\nu}(q) \COMMA
 \]
and includes hidden constraints on velocities $\dot{z} = \nu \Dq^T \xi \,\dot{q}$
implied by the constraints $(\pen{C})$ on position variables.
Replacing $\dot{z}$ and $z$ in $L_{\M{IMMP}}$ by their expressions as functions of $\dot{q}$ and $q$,
one obtains the Lagrangian associated with $\pen{H}$.
\end{proof}

The stochastically perturbed equations of motion of the Langevin type associated with \eqref{e:Himmp} define
the dynamics with implicit mass-matrix penalization.
\begin{Def}[IMMP]\label{d:IMMP} The implicit Langevin process associated with Hamiltonian $H_{\M{IMMP}}$
and constraints $(\pen{C})$ is defined by the following equations of motion
\begin{equation} \label{e:IMMP}\syst{
    &\dot{q} =  \MASST^{-1} p  &  \\
    &\dot{z} =  \MASSTZ^{-1} p_{z} &  \\
    &\dot{p} = -\Dq V(q) - \Dq V_{\M{fix},\nu}(q) -\gamma \dot{q} + \sigma \dot{W} - \Dq \xi \,\dot{\lambda} & \\
    &{\dot{p}_{z}} = -\gamma_{z} \dot{z} + \sigma_{z} \dot{W}_{z} +  \frac{\dot{\lambda}}{\nu} & \\
    &\xi(q) = \frac{z}{\nu}\COMMA & \TAG{\pen{C}}
  }
\end{equation}
  The process $\dot{W}$ (resp. $\dot{W}_{z}$ ) is a standard multi-dimensional white noise, $\gamma$ (resp. $\gamma_{z}$)
  a $d\times d$ (resp. $\nc\times \nc$) non-negative symmetric dissipation matrix,
  $\sigma$ (resp. $\sigma_{z}$) is the fluctuation matrix satisfying $\sigma \sigma^{T} =\frac{2}{\beta} \gamma$
  (resp. $\sigma_{z} \sigma_{z}^{T} =\frac{2}{\beta} \gamma_{z}$).
  The processes $\lambda \in \mathbb{R}^{\nc}$ are Lagrange multipliers associated with the constraints $(\pen{C})$ and
  adapted with the white noise.
\end{Def}

This process is naturally equivalent to the explicit mass-penalized Langevin process in $\mathbb{R}^{d}\times\R^d$
associated with $\pen{\MASST}$. Moreover, when the penalization vanishes ($\nu \to 0$), the evolution law of
the process $\{p_t,q_t\}_{t\geq 0}$ or $\{({p_{\nu}})_t,q_{t}\}_{t\geq 0}$ converges towards the original dynamics.
\begin{Pro}\label{p:IMMPeq} The stochastic process with constraints \eqref{e:IMMP} is well-posed and equivalent
to the Langevin diffusion in $\mathbb{R}^{d}\times\R^d$ (see Definition~\ref{d:langevin}),
with the mass-penalized Hamiltonian $\pen{H}$ \eqref{e:Hnu}, and the dissipation matrix given by
\[
\pen{\gamma}(q) = \gamma + \nu^{2} \Dq\xi \,\gamma_{z} \nabla^{T}_{q}\xi \PERIOD %
\]
Furthermore, the process is reversible and ergodic with respect to the canonical distribution \eqref{e:penboltz} (with marginal in position variables given by the original potential, i.e., up to the normalization,
$\EXP{ - \beta V(q) } dq$) .
\end{Pro}
\begin{proof}
Imposing the  constraints implies $ \nabla_q^T \xi\, \MASST^{-1} p = \frac{1}{\nu} \MASSTZ^{-1} p_{z}$.
Thus by the definition of $\pen{p}$ we have
\[
 \pen{p} = p + \nu \Dq\xi \, p_z \PERIOD
\]
 Since the position process $\{q_{t}\}_{t\geq 0}$ is of finite variation, a short computation shows that for each coordinate $i = 1,..,d$
\be\label{e:step}
\pen{\dot{p}}^i = {\dot{p}^i} + \nu (\partial_{q_i})^T \xi  \, \dot{p}_{z} + \nu^{2} \dot{q}^T   \Dq \! (\partial_{q_i}\xi) \, p_z \PERIOD
\ee
Furthermore,
\[
  -\partial_{q_i} \pare{\frac{1}{2}  \pen{p}^T \pen{\MASST}^{-1} \pen{p}}= \partial_{q_i}\pare{\frac{1}{2} \dot{q}^T\pen{\MASST} \dot{q}}
  = \nu^2 \dot{q}^T  \Dq(\!\partial_{q_i} \xi) \MASSTZ\nabla_q^T \xi\, \dot{q} \COMMA
\]
and thus
\[
{\pen{\dot{p}}} = \dot{p} + \nu \Dq \xi \,{\dot{p}_{z}} - \Dq \pare{\frac{1}{2}  \pen{p}^T \pen{\MASST}^{-1} \pen{p}} \PERIOD
\]
Substituting the expressions for  $\dot{p}$ and ${\dot{p}_{z}}$ from \eqref{e:IMMP} into \eqref{e:step} we obtain
\begin{equation}\label{e:IMMPeq}
{\pen{\dot p}}  =  -\frac{1}{2} \pen{p}^{T} \Dq \pen{\MASST}^{-1} \pen{p} - \Dq V(q)  - \Dq V_{\M{fix},\nu}(q)
\, -\gamma \dot{q} - \nu \Dq \xi\, \gamma_{z} \dot{z} + \sigma \dot{W} + \nu \Dq \xi\, \sigma_{z} \dot{W}_{z} \COMMA
\end{equation}
which yields the result.
\end{proof}
\begin{Pro}[Small penalty]\label{p:vanishpen}
When $\nu \to 0$ the evolution law of the processes $\{p_t,q_t\}_{t\geq 0}$ or
$\{(\pen{p})_t,q_t\}_{t\geq 0}$ defined by the implicit equations \eqref{e:IMMP} converges
(in the sense of probability distributions on continuous paths endowed with the uniform convergence)
towards the process solving the original Langevin dynamics \eqref{e:langevin}.
\end{Pro}
 \begin{proof}
The stochastic differential equation defined by $\dot{q} = \pen{M}^{-1} \pen{p}$
and \eqref{e:IMMPeq} has smooth coefficients which depend on $\nu$ in a continuous fashion ($\nu \mapsto \pen{M}$
and $\nu  \mapsto V_{\M{fix},\nu}$ are continuous).
Standard results on weak convergence (\cite{EthKur86}) of stochastic processes imply the result as stated.
\end{proof}

When the mass penalty tends to infinity, the IMMP process converges to a constrained process on the manifold $\SMANZ{z_{t=0}}=\set{ q \SEP \xi(q)=z_{t=0}}$ (surface measures are described in Section~\ref{s:measures}).
\bpro[Large penalty]\label{p:largepen}
        Consider a family of initial conditions indexed by $\nu$ and satisfying
	\[
	\sup_{\nu} \abs{ \nu\,(\xi(q_{t=0}) - z_{t=0}) } < +\infty \COMMA
	\]
        and assume that the Gram matrix $G$ is invertible in a neighborhood of $ \SMANZ{z_{t=0}}$.
	Then when $\nu \to +\infty$ the IMMP Langevin stochastic process \eqref{e:IMMP} converges weakly
        towards the decoupled limiting processes with constraints
	\begin{equation}\label{e:largepen}
        \syst{
    	  &\dot{q} =  \MASST^{-1} p \COMMA &  \\
    	  &\dot{p} = -\Dq V  - \Dq V_{\M{fix}} -\gamma \dot{q} + \sigma \dot{W} - \Dq\xi \dot{\lambda}\COMMA  & \\
    	  &\xi(q)  = z_{t=0}    \COMMA                                        &  \TAGG{C} \\
    	  &\dot{z} =  \MASSTZ^{-1} p_{z} \COMMA &  \\
    	  &{\dot{p}_{z}} = -\gamma_{z} \dot{z} + \sigma_{z} \dot{W}_{z} \PERIOD &
	}
        \end{equation}
	where $\{\lambda_t\}_{t\geq 0}$ are adapted stochastic processes defining
        the Lagrange multipliers associated with the constraints $(C)$.

	Furthermore, the process $\{q_t,p_t\}_{t\geq 0}$ defines an effective dynamics with constraints
        (see also Definition~\ref{d:consteff}) on the sub-manifold $\SMANZ{z_{t=0}}$. It is reversible with respect to its stationary canonical distribution given, up to the normalization, by
        \[
	\EXP{-\beta (H(p,q)+V_{\M{fix}}(q))} \sigma_{T^{*}\SMANZ{z_{t=0}}}(dp\,dq)
	\]
	with the $q$-marginal $\EXP{- \beta V(q)}\,\delta_{\xi(q) = 0}(dq)$.
	When $\gamma$ and $\gamma_{z}$ are strictly positive definite the process is ergodic.
\end{Pro}
\begin{proof}
By a simple translation, it is sufficient to show the proposition for $z_{t=0} = 0$. Satisfying the constraint $(\pen{C})$ in \eqref{e:IMMP} implies a hidden constraint in the momentum space,
$\nabla_q \xi \MASST^{-1} p$ $= \frac{1}{\nu} \MASSTZ^{-1} p_{z}$. Differentiating this expression with respect to time
and replacing the result in \eqref{e:IMMP} yields an explicit formula for the Lagrange multipliers
\be\label{e:lagecd}
  \dot{\lambda} = (G+\frac{1}{\nu^2}\MASSTZ^{-1})^{-1}\left[\HESS(\xi) \pare{\MASST^{-1} p,\MASST^{-1} p} +
                                                             \Dq \xi \MASST^{-1} f_{q}
                                                           - \frac{1}{\nu}\MASSTZ^{-1} f_{z}\right] \COMMA
\ee
with forces $(f_q,f_z)$
\begin{eqnarray*}
f_{q} &=&  -\Dq V - \Dq V_{\M{fix},\nu}  -\gamma \MASST^{-1} p + \sigma \dot{W}\COMMA \\
f_z &=&   - \gamma_{z} \MASSTZ^{-1} p_{z} + \sigma_{z} \dot{W}_{z}\COMMA
\end{eqnarray*}
and the Hessian $\HESS(\xi)$ of the mapping $\xi$ acting on the velocities $\MASST^{-1}p$.
This calculation shows that \eqref{e:IMMP} is in fact a standard stochastic differential equation with smooth coefficients,
and thus has a unique strong solution. The coefficients of these stochastic differential equations are
continuous in the limit $\tfrac{1}{\nu} \to 0$, at least in a $\delta$-neighborhood of $\SMAN$ in
which $G$ is invertible.
The formally computed  limiting process is given by \eqref{e:largepen} with the Lagrange multipliers solving
 \[
 \dot{\lambda} = G^{-1}\pare{\HESS(\xi) \pare{\MASST^{-1} p,\MASST^{-1} p} + \Dq^T \xi \, \MASST^{-1} f_{q}} \PERIOD
 \]
By construction, this limiting process satisfies the constraint $\xi(q)=0$.
Its coefficients are Lipschitz and the process is well-posed.
As a result of those properties, the rigorous proof of weak convergence follows
classical arguments, see \cite{EthKur86}, that are divided into three steps
\begin{description}
\item[{\rm (i)}] We truncate the process \eqref{e:IMMP} to a compact neighborhood of $\mathcal{M}_{0}$.
\item[{\rm (ii)}] The continuity of the Markov generator with respect to $\tfrac{1}{\nu}$ implies tightness for the associated $\tfrac{1}{\nu}$-sequences of truncated processes with the limit being uniquely defined by \eqref{e:largepen}.
\item[{\rm (iii)}] The limiting process remains on  $\SMAN$, which implies weak convergence of the sequence without truncation.
\end{description}
The process \eqref{e:largepen} is thus a Langevin process with constraints, exhibiting reversibility properties
with respect to the associated Boltzmann canonical measure and is ergodic when $\gamma$ is strictly positive definite
(see the summary in Appendix~\ref{s:langevinapp}).
Note that the $q$-marginal is geometrically corrected by the Fixman potential term.
\end{proof}

We conclude this section by discussing some consequences for numerical computations.
\begin{Rem}
\begin{description}
\item[{\rm (i)}] Proposition~\ref{p:vanishpen}  and ~\ref{p:largepen} imply that IMMP scheme is a tunable
                 interpolation between the exact dynamics \eqref{e:langevin}, and the dynamics with constraints \eqref{e:largepen}.
\item[{\rm(ii)}] In Section~\ref{s:highoscill} where a stiffness paramater, corresponding to a large
                 penalty $\nu$, is introduced
                 for slow/fast systems leads to a Markovian effective motion constrained on
                 the slow manifold $\SMAN$, usually seen as a heuristically consistent
                 Markovian approximation of the exact dynamics
                 (see \cite{Rei00}).
\end{description}
\end{Rem}

\section{Numerical integration}\label{s:numerinteg}
The key ingredient for achieving efficient numerical simulation is to use an integrator
that enforces the constraints associated with the implicit formulation of the mass penalized dynamics \eqref{e:IMMP}. The implicit structure of \eqref{e:IMMP} leads to numerical schemes that are potentially asymptotically stable in stiff cases (Section~\ref{s:highoscill}). On the other hand, when the penalization $\nu$ vanishes with the time-step
the scheme becomes consistent with respect to
the original exact dynamics \eqref{e:langevin}. One may then consider the mass-penalization introduced here as
a special method of
pre-conditioning for a stiff ODE system with an ``implicit'', in the time evolution sense, structure.
Here, the ``implicit'' structure amounts to solving the imposed constraints $\xi(q) = \frac{z}{\nu}$ in \eqref{e:IMMP}.

It lies outside the scope of this paper to make a review of standard numerical methods for constrained mechanical systems,
we refer to \cite{HaiLub02} as a classical textbook, and to the series (\cite{GunBer77, RycCic77, CicVan06, CicKapVan05, CicLelVan06})
as a sample of works on practical developments of numerical methods.
The IMMP method is presented  with the classical leapfrog/Verlet scheme that enforces constraints,
usually called RATTLE in \eqref{e:immpscheme}. It can be implemented by a simple modification of standard schemes constraining fDOFs. The scheme is second order, reversible and symplectic. This choice is largely
a presentation matter,
for practical purposes one can refer to one's favorite numerical integrator for Hamiltonian systems
with or without stochastic perturbations.

As an option, one can also add a Metropolis acceptance/rejection corrector at each time step. If the underlying integrator is reversible and preserves the phase-space measure, this extension leads to a scheme which {\it exactly} preserves canonical distributions. Such algorithms are usually referred to as Hybrid Monte Carlo (HMC). Generically  HMC are sampling algorithms which use the underlying dynamics of
the system to generate moves in the configuration space, which are accepted or rejected according to a Metropolis rule.  While this method may allow for larger time steps in the integrator it may also lead to higher computational costs when the
acceptance ratio becomes low. The acceptance ratio often decreases when the dimension of the system increases
(\cite{IzaHam04}). However, many improvements and modifications (\cite{Creutz89, Horowitz91, Beccaria94, Kennedy98})
of the HMC algorithm  have been developed since its
introduction in \cite{Duane85, DuaneKennedy87}. While the HMC methods were
initially used in simulations applied to quantum statistical field theories, they have been also employed to wide range
of simulations in macromolecular systems, e.g., \cite{ClaBakStiBra94,SchFis99}.

We implement numerical discretization of the Langevin process with constraints \eqref{e:IMMP} obtained
by splitting the Hamiltonian part and the Gaussian fluctuation/dissipation perturbation.
The idea of using a Hybrid Monte Carlo step in discretizing a Langevin process goes back to \cite{Horowitz91}.
Note that by using a Metropolis acceptance/rejection rule (HMC), or simply by weighting statistical averages,
the forces associated with the Fixman corrector need not be computed when one is only interested in sampling.

We recall that we consider the IMMP dynamics \eqref{e:IMMP}, which consists of a Langevin process
with constraints defined by the following elements
\begin{enumerate}
\item the Hamiltonian $H_{\M{IMMP}}$ defined in \eqref{e:Himmp},
\item the dissipation matrix
$\begin{pmatrix}
\gamma & 0 \\
0  & \gamma_{z}
\end{pmatrix}$,
\item the inverse temperature $\beta$.
\end{enumerate}
\begin{Rem}
A useful variant, when using HMC strategies, consists in modifying the Hamiltonian $H_{\M{IMMP}}$ \eqref{e:Himmp} in the integrator (importance sampling) by neglecting the Fixman corrector
\begin{equation}\label{e:Hnum}\syst{
&\tilde{H}_{\M{IMMP}}(p,p_{z},q,z) =  \frac{1}{2} p^{T} \MASST^{-1} p  + \frac{1}{2} p_{z}^{T} \MASSTZ^{-1} p_{z} + \tilde{V}(q)  & \\
&\xi(q) = \frac{z}{\nu}         & \TAG{\pen{C}},
}\end{equation}
where $\tilde{V}$ is a potential that may be chosen arbitrarily (typically $\tilde{V} = V$). Indeed, only the underlying phase space structure, which does not depend on $\tilde{V}$, is necessary in HMC methods. The correct potential $V + V_{\M{fix},\nu}$ has only to be used in the Metropolis step, or by weighting ensemble averages, to ensure exact canonical sampling. Thus potentially costly evaluations of the gradient of the Fixman corrector $V_{\M{fix},\nu}$ can be avoided.
\end{Rem}
\begin{scheme}[Leapfrog/Verlet algorithm for Langevin IMMP \eqref{e:IMMP}]\label{d:scheme}

\bi

\item[Step 1:] Integrate the Hamiltonian part with:
\be \label{e:immpscheme}
                \arr{
		&&\syst{
		&p_{n+1/2} = p_{n} -  \frac{\dt}{2} (\Dq V+\Dq V_{\M{fix},\nu})(q_{n}) - \Dq \xi(q_{n})  \lambda_{n+1/2}&\\
		&p_{n+1/2}^z = p_{n}^z +\frac{1}{\nu}\lambda_{n+1/2} &
		}\\
		&&\syst{
		&q _{n+1} =  q_{n}+ \dt \MASST^{-1} p_{n+1/2} &\\
		& z_{n+1} =  z_{n}+ \dt \MASSTZ^{-1} p_{n+1/2}^z &
		}\\
		&&\hspace{.5cm} \xi(q_{n+1}) = \frac{z_{n+1}}{\nu} \qquad \hspace{3cm} (C_{1/2}) \\
		&&\syst{
		&p_{n+1}   = p_{n+1/2} -  \frac{\dt}{2} (\Dq V+\Dq V_{\M{fix},\nu})(q_{n+1})  - \Dq \xi(q_{n+1})  \lambda_{n+1}&\\
		&p_{n+1}^z = p_{n+1}^z +\frac{1}{\nu}\lambda_{n+1}&
		}\\
		&&\hspace{.5cm} \Dq^T \xi(q_{n+1})\,\MASST^{-1}p_{n+1} = \frac{1}{\nu} \MASSTZ^{-1}p^z_{n+1}\;\;\;\;\;
                 (C_{1})   \PERIOD \\
		}
\ee
\item[Step 2:] Integrate the Gaussian fluctuation/dissipation part with a mid-point Euler scheme with constraints (see Appendix~\ref{s:exactflucdiss}).
\ei
\end{scheme}
To obtain exact sampling by correcting the time step errors, one needs to introduce a Metropolis correction.
\begin{scheme}[HMC]\label{d:schemehmc}
\bi

\item[Step 1:] Compute $(q_{n+1},z_{n+1},p_{n+1},p^{z}_{n+1})$ with \eqref{e:immpscheme}, and set $ \Delta H_{n+1} = H_{\M{IMMP}}(q_{n+1},z_{n+1},p_{n+1},p^{z}_{n+1}) - H_{\M{IMMP}}(q_{n},z_{n},p_{n},p^{z}_{n})$.

\item[Step 2:] Accept the step with the probability $\min(1,{e \rm }^{-\beta \Delta H_{n+1} }$, otherwise reverse impulses and set
\[
(q_{n+1},z_{n+1},p_{n+1},p^{z}_{n+1})= (q_{n},z_{n},-p_{n},-p^{z}_{n}) \PERIOD
\]

\item[Step 3:] Integrate the Gaussian fluctuation/dissipation part with a mid-point Euler scheme (for details see Appendix~\ref{s:exactflucdiss}).
\ei
\end{scheme}
This numerically constructed Markov chain preserves the canonical distribution.
\bpro[Exact sampling]
Assume that \eqref{e:immpscheme} is defined globally. The numerical discretization of \eqref{e:IMMP}
described in Scheme~\ref{d:schemehmc} generates a Markov chain that leaves the canonical
distribution \eqref{e:penboltz} invariant, with marginal distribution in position variables is the original distribution $\EXP{ - \beta V(q) } dq$, which is independent of the mass-penalization $\nu$.
\epro
\begin{proof}
The statement follows from reversibility and measure preserving properties of Verlet schemes (see \cite{HaiLub02}),  from the Hybrid Monte Carlo rule (\cite{DuaneKennedy87}),
and from the construction of the mass-penalized Hamiltonian (Proposition~\ref{p:IMMPeq}).
\end{proof}

One can now construct consistent schemes by letting the penalty $\nu =\nusc \dt^{k}$ go to zero with
the time-step, for some $k >0$.
Indeed, Proposition~\ref{p:vanishpen} shows that the mass-penalized dynamics \eqref{e:IMMP} converges towards the exact original dynamics for $\nu=0$ at order $\nu^2$. Consequently most of the usual numerical schemes will be consistent at  their own approximation order, but bounded above by $2k$. Neglecting the order of the fluctuation/dissipation part (Step 2 of Scheme~\ref{d:scheme}), we deduce the following convergence property. The order of convergence refers to the maximal integer $k$ such that the convergence of trajectories
with respect to the uniform norm occurs at the rate $\BIGO(\dt^{k})$.
\bpro[Time-step consistency]
Assume that the numerical flow \eqref{e:immpscheme} is defined globally, and that $\nu =  \nusc \dt^{k}$. Then the IMMP numerical scheme \eqref{e:immpscheme} is consistent of the order $\min(2,2k)$ with respect to the original exact limiting process \eqref{e:langevin} without thermostat ($\gamma=0$, $\sigma=0$).
\epro
\begin{proof}
We consider the scheme \eqref{e:immpscheme}	in the variables $(p,q,\nu p_z,\nu z)$.
Then the force field only depends on $q$ and the global mass-matrix is smooth with respect to $\nu^2$.
The sub-manifold defining the constraints ($\nu^2 \xi(q) = \nu z $) is also smooth with respect to $\nu^2$.
By the implicit function theorem we have that locally the RATTLE scheme is the standard leapfrog scheme (see \cite{HaiLub02,LeiRei05}), and the local mapping depends smoothly on $\nu^2$.
Therefore the standard calculation of the order of the leapfrog scheme holds uniformly with respect to $\nu$, \cite{HaiLub02}. Now, the mass-penalized dynamics \eqref{e:IMMP} is a differential equation (deterministic here) whose coefficients differ from
the original
process by a smooth perturbation of order $\nu^2$. The result thus follows from applying a simple Gronwall argument.
\end{proof}

Note also that taking $\nu \to +\infty$ in \eqref{e:immpscheme} gives a numerical scheme with constraints consistent with the limiting constrained dynamics of Proposition \ref{p:largepen}. For slow/fast systems (see Section~\ref{s:infstiff}), large penalty will lead to asymptotic stability in the stiffness limit.

We conclude this section with a practical recipe for tuning the mass matrix penalty. This can be done for instance by computing the time averaged Metropolis acceptance ratio $r= \EXPECT{{e \rm }^{-\beta (\Delta H_{n})}}$
in Scheme~\ref{d:schemehmc}, which gives a precise quantification of the time-step error.
Then increasing the penalty $\nu$ can save computational time as long as it leads to an increase 
of the average ratio $r$. Indeed, this means that the selected fDOFS are limiting the time-step stability region. Prescribing a time average ratio $r$ (for instance $0.90$), a maximal time-step $\dt_{max}$ 
associated with the largest penalty $\nu_{max}$ that is able to improve stability can be obtained in this way.
Finally, one can set, for example, $\nu = \tfrac{\nu_{max}}{\dt_{max}} \dt$ in Scheme~\ref{e:immpscheme} 
to obtain an order $2$ convergent scheme with increased stability region.

\section{A high dimensional numerical test case}\label{s:numerexperiment}
The behavior of numerical methods for Hamiltonian systems becomes of particular interest when simulating
high-dimensional systems.
Such systems in general exhibit many different time-scales and one seeks to go beyond stability constraints
implied by the fastest microscopic degrees of freedom. As mentioned in the introduction, globally implicit methods are usually too costly, while efficient splitting methods using minimal implicitness represent an active research area.
In order to demonstrate applicability of the IMMP method
we consider non-convex repulsive interactions between particles in the gas phase which makes
implicit methods intractable.

On the other hand, direct constraints on microscopic fast degrees of freedom are commonly used in practice and can introduce a potentially strong bias in the macroscopic behavior of the whole system.
In the case presented here, a particle chain with pairwise interactions, one may be interested in some
macroscopic observable, like the total length of the chain.
Direct constraints would then lead to a \emph{totally rigid chain}, losing completely the evolution of
the latter macroscopic variable.

We numerically study the IMMP strategy described in Section~\ref{s:IMMP} and Section~\ref{s:numerinteg}, and  the latter is systematically compared with the exact dynamics numerically integrated by a simple Leapfrog/Verlet scheme with a thermostat.
Langevin thermostats are used, with a Metropolis/HMC step for the dynamics integrator as described by Scheme \ref{d:schemehmc}. An appropriate scaling of the chain of size $N$ at the mass transport level with a large mass-penalty of order $\nu :=\nusc N $ are chosen. Numerical results show that IMMP induces a gain in time stepping
of order $N$, while the macroscopic dynamics (in particular convergence to equilibrium), given by a formal non-linear stochastic partial differential equation,
remains of order $1$. Macroscopic dynamics of the IMMP method also stays close to the exact Verlet integration when re-scaled penalty $\nusc$ is small.

\subsection{The particle chain model}
The model we consider consists of a chain of particles which interact through a non-convex repulsive
pairwise potential of the form
\begin{equation}\label{e:vint}
 \syst{
& v_{\INTP}(r) = f(r)\COMMA    \quad r\leq r_{\INTP}\COMMA   & \\
& v_{\INTP}(r) = 0\COMMA       \quad\quad r >  r_{\INTP} \COMMA   &
}
\end{equation}
where $ f(r_{\INTP}) = f'(r_{\INTP})=0$ and $\ds \lim_{r \to -\infty} f(r) = +\infty$. Each particle is also individually submitted to a macroscopic confining exterior potential.
After converting to the non-dimensional form the typical quantities involved in the model enable us to write
a scaling at the mass-transport level where the dynamics of the chain is described by the Hamiltonian
\be\label{e:HN}
H_{N}(q,p) = \frac{1}{2} p^{T}  p + \sum_{i=1}^{N-1}  v_{\INTP}(\dg_{i} q ) + \sum_{i=1}^{N} v_{\EXTP}(q_{i})  \COMMA
\ee
and by a coupling with an exterior thermal bath at the re-scaled inverse temperature
\begin{equation}\label{e:beta}
 \beta_N :=\beta N^{-1} \PERIOD
\end{equation}
In the expression \eqref{e:HN} the functions $r \in\R \mapsto v_{\INTP}(r)\in\R$ and $q \mapsto v_{\EXTP}(q)\in\R$
are the smooth interaction potential
and the exterior potential rrespectively
The linear operator $\dg:\mathbb{R}^{N} \to \R^{N-1}$ having the components
 \[
 \dg_{i} q = \frac{q_{i+1}-q_{i}}{1/N}\COMMA\;\;\;i=1,\dots,N-1\COMMA
 \]
represents the discrete gradient associated to the chain with the Neumann boundary conditions.
Its transpose operator is denoted  $(\dg)^T:\mathbb{R}^{N-1} \to \R^{N}$.
The particles are represented by their re-scaled positions $q=(q_1,...,q_N)$, so that the typical position
and deviation of  $q$ is formally of order $1$ with rrespectto $N$.
This can be seen by considering particles in the chain as indexed by $x= \tfrac{i}{N} \in [0,1]$.
The thermodynamic limit $N\to\infty$ of the Hamiltonian \eqref{e:HN} then %
converges to
 \begin{equation}\label{e:Hpde}
 \lim_{N \to \infty } \frac{1}{N}H_{N}(q,p) = \int_{0}^{1} \frac{1}{2} p(x)^{2} \,dx + \int_{0}^{1} v_{\INTP} \left(\frac{dq}{dx}\right) \,dx +
         \int_{0}^{1} v_{\EXTP}(q(x)) \,dx \COMMA
 \end{equation}
when $x\mapsto q(x) $ is a smooth function on $[0,1]$ satisfying the Neumann boundary conditions
$q'(0)=q'(1)=0$.
We choose to work with the scaling given by \eqref{e:HN}-\eqref{e:Hpde} in order to prescribe the macroscopic timescale of the chain profile at order $1$ with respect to $N$.
Using the matrix notation, we write the stochastically perturbed equations of motion
\[
\arr{
\dot{q} &=&  p \\
\dot{p} &=& -(\dg)^T v'_{\INTP}( \dg q ) - v'_{\EXTP}(q)- \gamma p + \sqrt{N}\sigma\dot{W}
}
\]
with the fluctuation/dissipation identity $\sigma \sigma^{T} = 2\beta^{-1}\gamma$.

Next we turn to the mass-penalized Hamiltonian model for this particle chain.
The number of independent fast degrees of freedom
that will be penalized is equal to $N-1$.  The penalized degrees of freedom 
are given by the inter-particle distances
$\xi_{i}(q) = z_{i} = q_{i}-q_{i-1}$.
The mass-penalization intensity has to scale formally as
\[
\nu:=\nu_{N} =\nusc N
\]
to give a discrete operator meaning to the perturbed mass-matrix in $\mathbb{R}^{N}$
\[
\penN{\MASST} = \ID + \nusc^{2}  (\dg)^{T} \MASSTZ \dg \PERIOD
\]
When $\MASSTZ=\ID$ the matrix perturbation is given by the discrete Laplacian
\[
\dl := - (\dg)^{T} \dg \PERIOD
\]
with the Neumann boundary conditions. Following our general construction we obtain the mass-penalized
Hamiltonian
\be\label{e:HNpen}
\penN{H}(\penN{p},q) = \frac{1}{2} \penN{p}^{T} \MASST^{-1}_{\nu_N} \penN{p}
                       + \sum_{i=1}^{N}  v_{\INTP}( \dg_i q )
                       + \sum_{i=1}^{N} v_{\EXTP}(q_{i}) \COMMA
\ee
and the stochastically perturbed equations of motion with $\MASSTZ =\ID$ are
\be \label{e:spde}
 \ddot{q} =  ( \ID -  \nusc^{2} \dl  )^{-1}  \left( -(\dg)^T v'_{\INTP}( \dg q )
               - v'_{\EXTP}(q) - \gamma \dot{q}
               +   \sigma \sqrt{N}\dot{W} \right)  \PERIOD
\ee
The form \VIZ{e:spde} of the equations of motion clearly demonstrates
the regularizing effect of the mass-matrix penalization.
Indeed, all the high frequencies generated by the microscopic forces
$-(\dg)^T v'_{\INTP}( \dg q )$ are  \emph{filtered} by a ``low-pass filter'' represented
by the regularizing operator $\pare{ \ID - \nusc^{2} \dl }^{-1}$.

\subsection{Numerical results}
In this section we demonstrate behavior of the IMMP method with the Hamiltonian \eqref{e:HNpen}
for different values of the re-scaled mass-matrix penalization intensity $\nusc$, the time step parameter $\dt$,
and the system size $N$.
In particular, we perform a systematic comparison with the exact dynamics \eqref{e:HN} integrated by
the Leapfrog/Verlet scheme.

The simulations were carried out with the re-scaled quadratic exterior potential
\[
v_\EXTP(q) = \left( \frac{q-0.5}{2.2}\right)^2 \COMMA
\]
and with a double-well repulsive re-scaled interaction potential
\[\syst{
& v_{\INTP}(r) = ( 50\pare{ (r - 0.1)^2 - 0.05^2})^2\COMMA      \quad r\leq 0.1 & \\
& v_{\INTP}(r) = 0\COMMA       \quad r >  0.1 \COMMA   &
}
\]
and the rescaled inverse temperature $\beta_N =10\,N^{-1}$.
We choose the penalizing mass-matrix to be the identity matrix $\MASSTZ = \ID$,
and the dissipation parameter $\gamma=0.1$.

\medskip
\noindent
{\sc Test I:} {\it Comparison of macroscopic dynamics.}

\begin{figure}[ht]
  \centerline{
  \includegraphics[width=8.5cm]{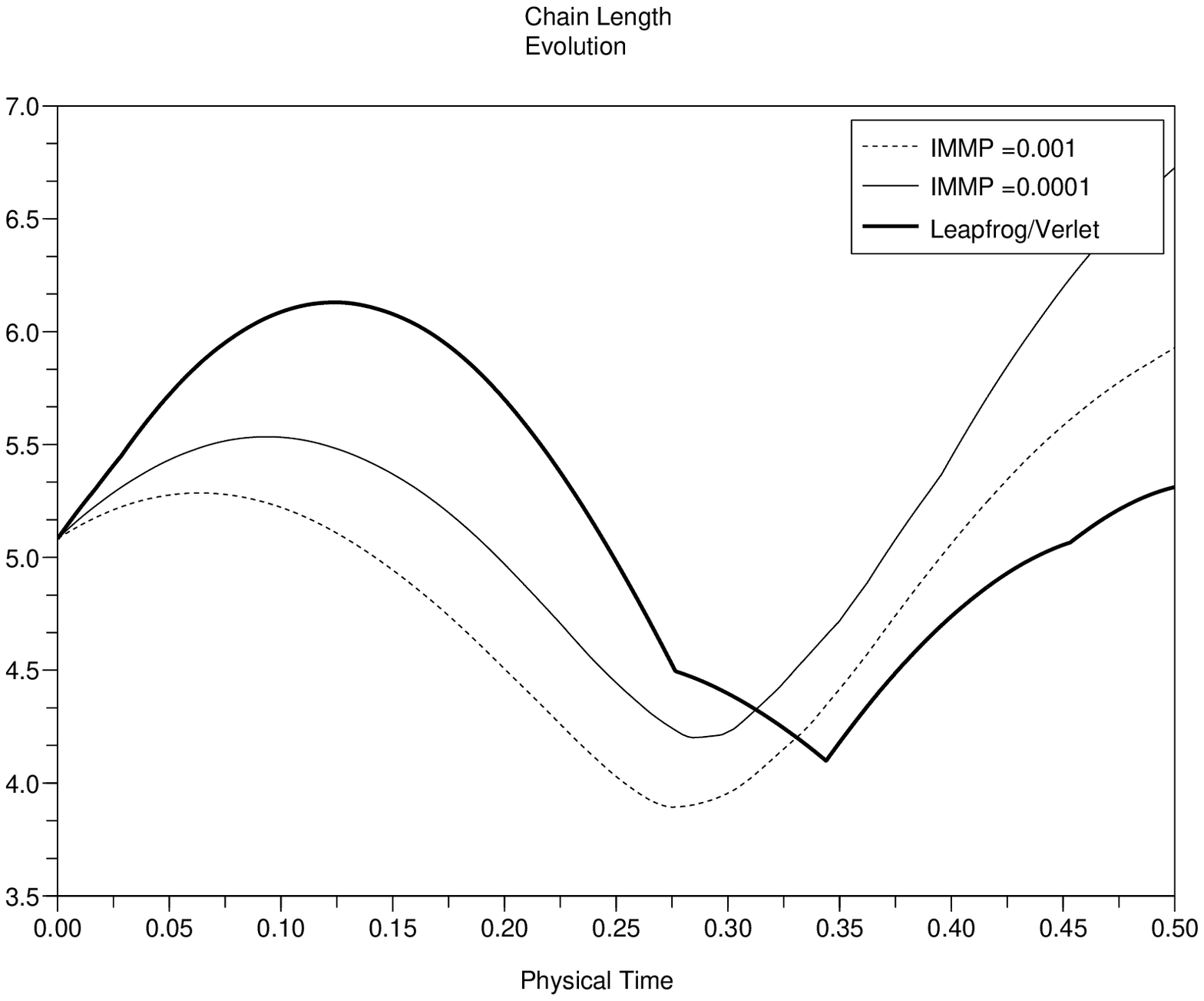} \hspace{-.3cm}
  \includegraphics[width=8.5cm]{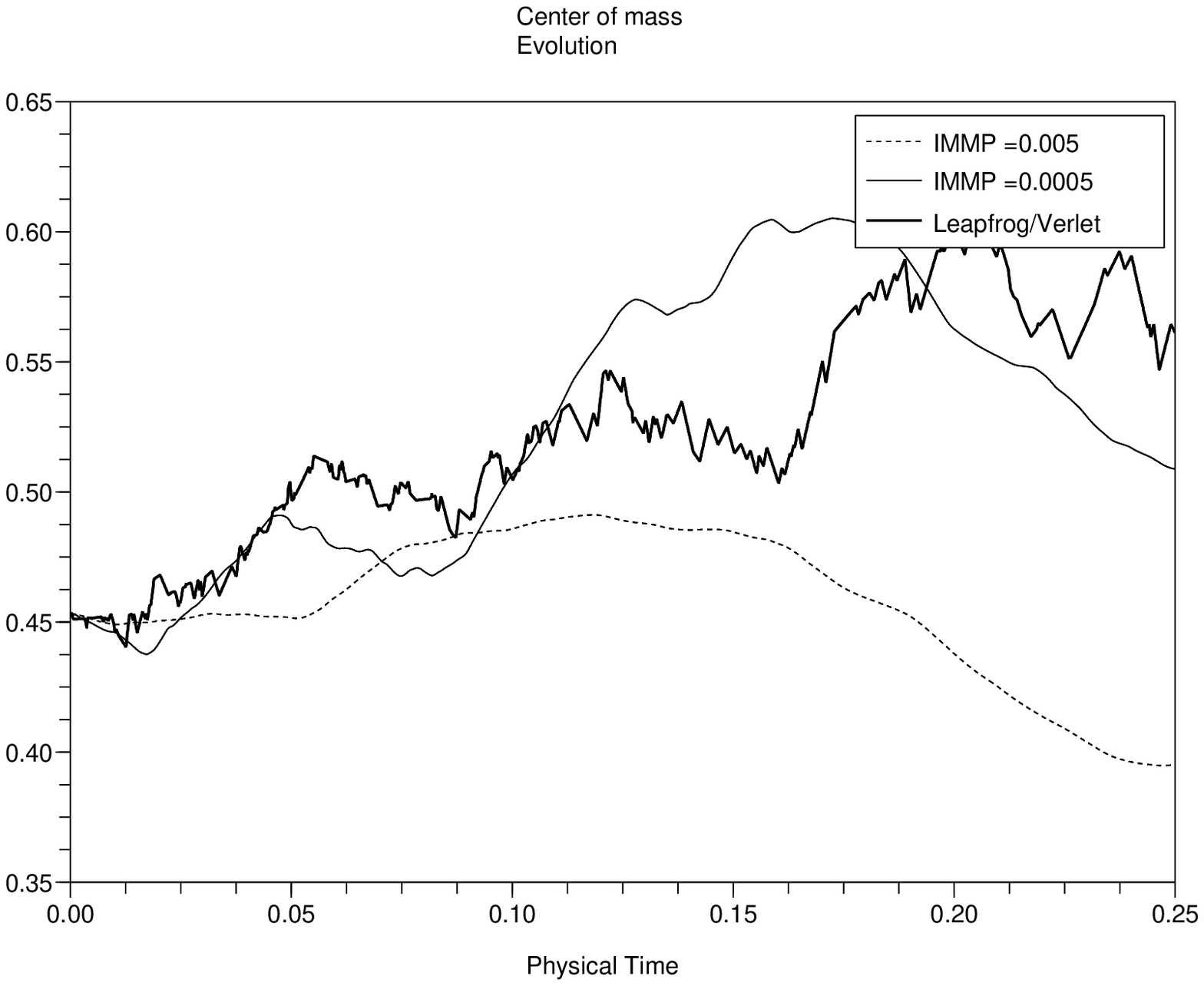}
}
  \caption{Short time evolution of the chain length and the center of mass for different values of
           the IMMP mass-penalization $\nusc^2 = \{10^{-3},10^{-4} \} $ and
           $\nusc^2 = \{5.10^{-3},5.10^{-4} \}$, $N=500$.}\label{f:pen}
\end{figure}

We demonstrate approximation of the macroscopic dynamics by the IMMP method simulating two observables:
(i) the length of the chain: $l = q_N - q_1$, and (ii) the center of mass: $c = q_{N/2}$. %
The dynamical behavior for short times of the chain length $l$ and the center of mass $c$ for different values of the mass-matrix penalization intensity $\nusc$ is depicted in Figure~\ref{f:pen}.
The same realization of the noise is used in each simulation, and no Metropolis rejection is used. One can observe the smoothing induced by the IMMP method, however, typical evolution timescales are not changed
by the IMMP method.

\begin{figure}[ht]
\centerline{
\includegraphics[width=8.5cm]{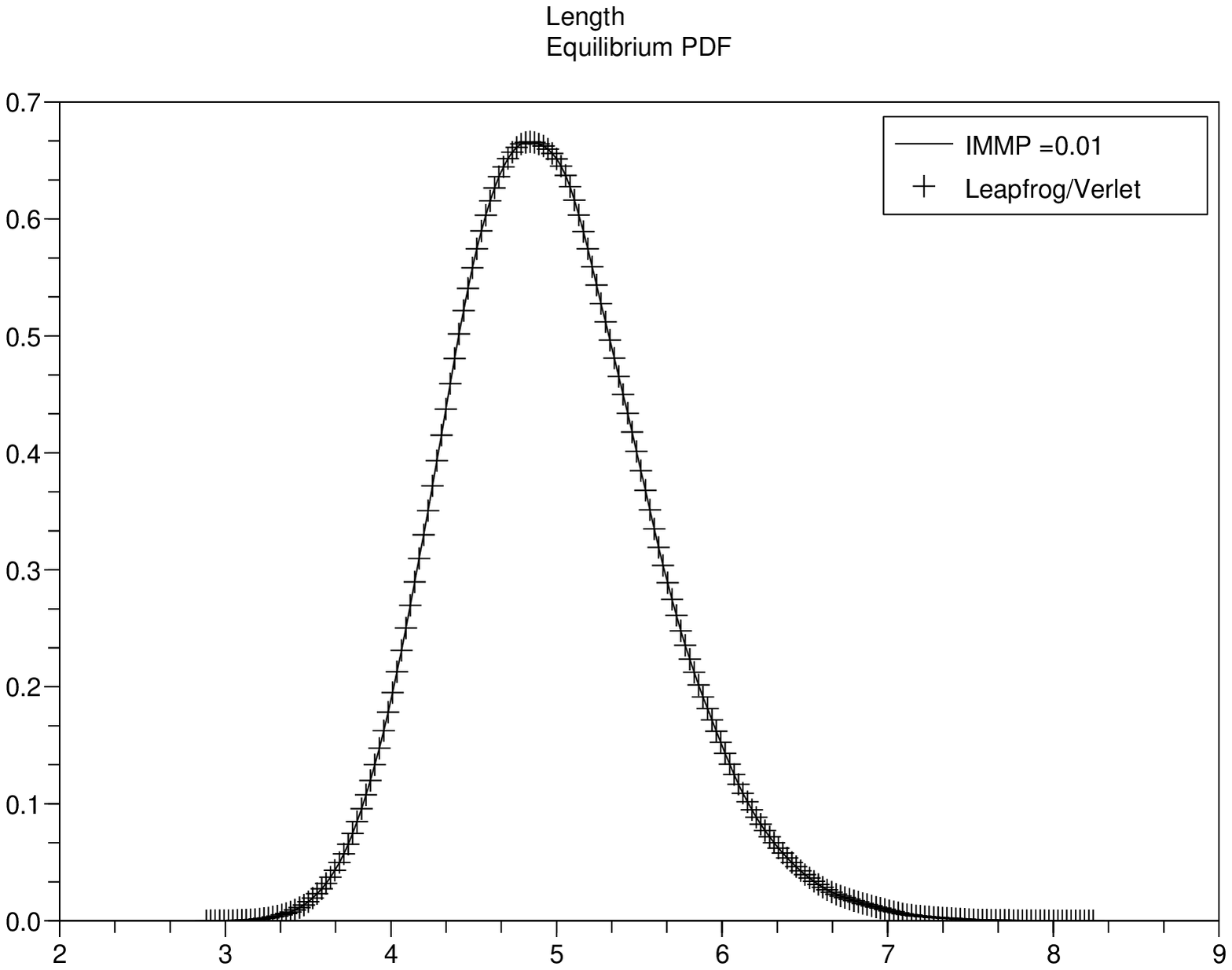}
\includegraphics[width=8.5cm]{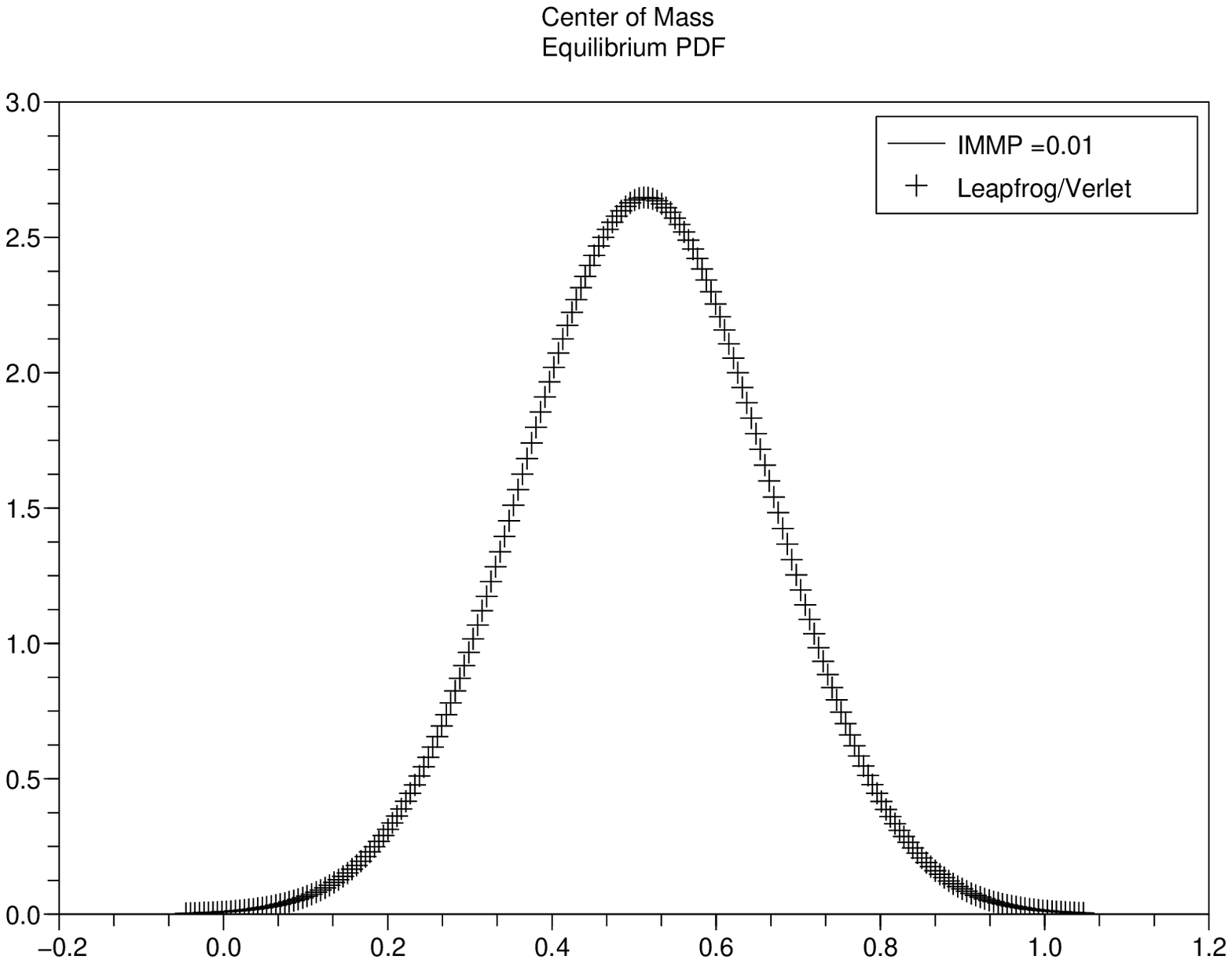}
}
\caption{Equilibrium PDF of the chain length and the center of mass. Simulated by the IMMP method with the
         penalty $\nusc^2 = 10^{-2}$.}\label{f:equilibrium}
\end{figure}

\begin{figure}[ht]
\includegraphics[width=8.0cm]{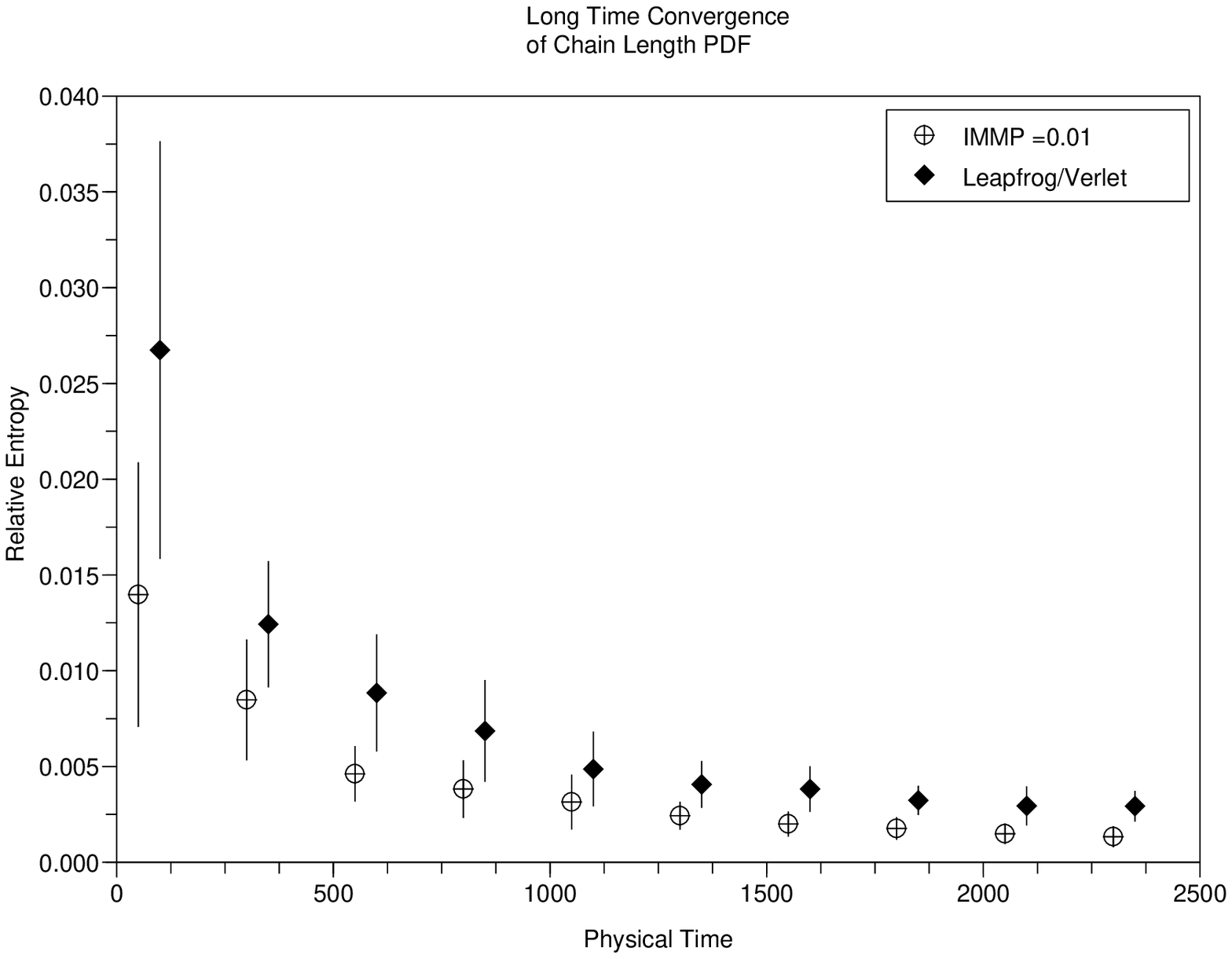}
\includegraphics[width=8.0cm]{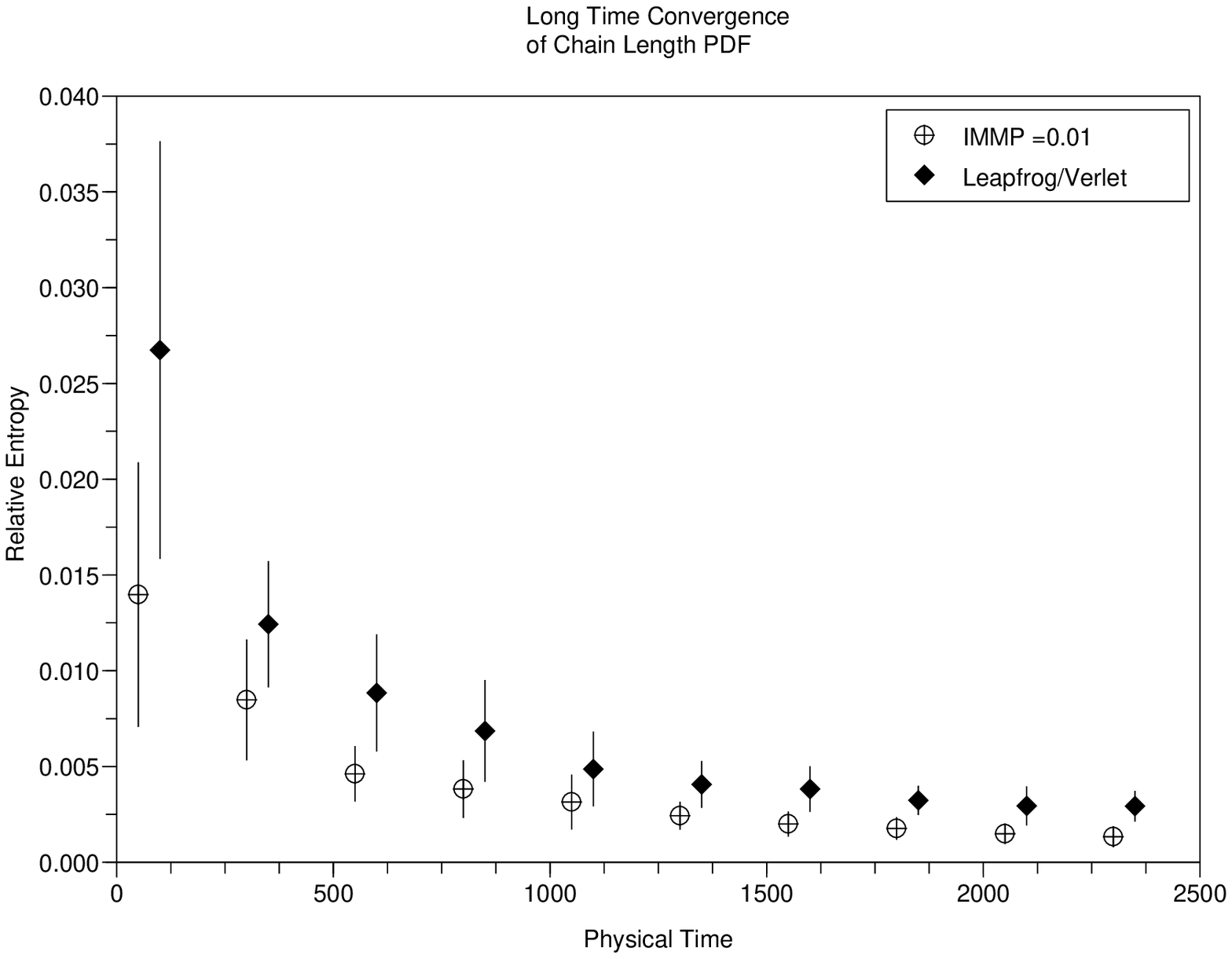}
\vspace{0cm}
\caption{Convergence of the relative entropy of the PDF of a trajectory of the chain length and the
         center of mass towards equilibrium. Simulated by the IMMP method with the
         penalty $\nusc^2 = 10^{-2}$. The error bar represents extreme quartiles as computed over $200$ realizations.}\label{f:conv}
\end{figure}

Next, the convergence to equilibrium of $l$ and $c$ are analyzed. The initial condition is taken out of equilibrium with zero potential energy and particles concentrated near the central position  $q = 0.5$. The equilibrium probability density functions (PDF) of $l$ and $c$ are computed for the IMMP case and the exact dynamics in Figures~\ref{f:equilibrium},
exact sampling ensuring identity of the two calculations. The PDFs are estimated using a Gaussian kernel estimator.
The evolution of the relative entropy between a trajectory sample and the equilibrium distribution is plotted in Figure~\ref{f:conv}. We observe that the convergence of the IMMP dynamics is slightly faster than the exact
dynamics with the Verlet integrator.
Convergence to equilibrium on shorter times can be also observed by inspecting the auto-correlations in
time series of $l(t)$ and $c(t)$ in Figure~\ref{f:conv}.
Convergence of the IMMP method is again faster. However, it is interesting to note that the mass ppenalization
introduces an inertial artifact that produces oscillations in the auto-correlations of the center of mass evolution.

The average transition time $\tau(\nusc, N)$ of the chain center of mass $c$ between $q=0.4$ and $q=0.6$ is summarized in
Table~\ref{t:avtime}\footnotemark[1] for
different values of the chain size $N$ and the mass penalization $\nusc$.
The results show that the average time $\tau(\nusc, N)$ of the macroscopic dynamics of $c$ scales
as $\BIGO(1)$ with respect to the system size, and is continuous with respect to the exact dynamics when  $\nusc \to 0$.

\begin{table}[h]
\begin{center}
\begin{tabular}{|c|ccccc|}
\hline
  & $\nusc^2 = 10^{-2}$ &  $\nusc^2 = 10^{-3}$  & $\nusc^2 = 10^{-4}$ & $\nusc^2 = 10^{-5}$  & Leapfrog/Verlet  \\
\hline
N = 100 & 1.6 (1) &  1.3 (1) & 1.2 (1)&  1.0 (1) &1.0 (1) \\
\hline
N = 300 & 1.5 (1) &  1.3 (1) & 1.2 (1)&  1.1 (1) &1.0 (1) \\
\hline
N = 500 & 1.4 (1) &  1.2 (1) & 1.2 (1)&  1.1 (1) &1.0 (1) \\
\hline
\end{tabular}
   \caption{\label{t:avtime} Average transition time of the chain center of mass from the position $q_{N/2}=0.4$
            to the position $q_{N/2}=0.6$, simulated by the IMMP method.}
\end{center}
\end{table}

\begin{figure}[ht]
\includegraphics[width=8.0cm]{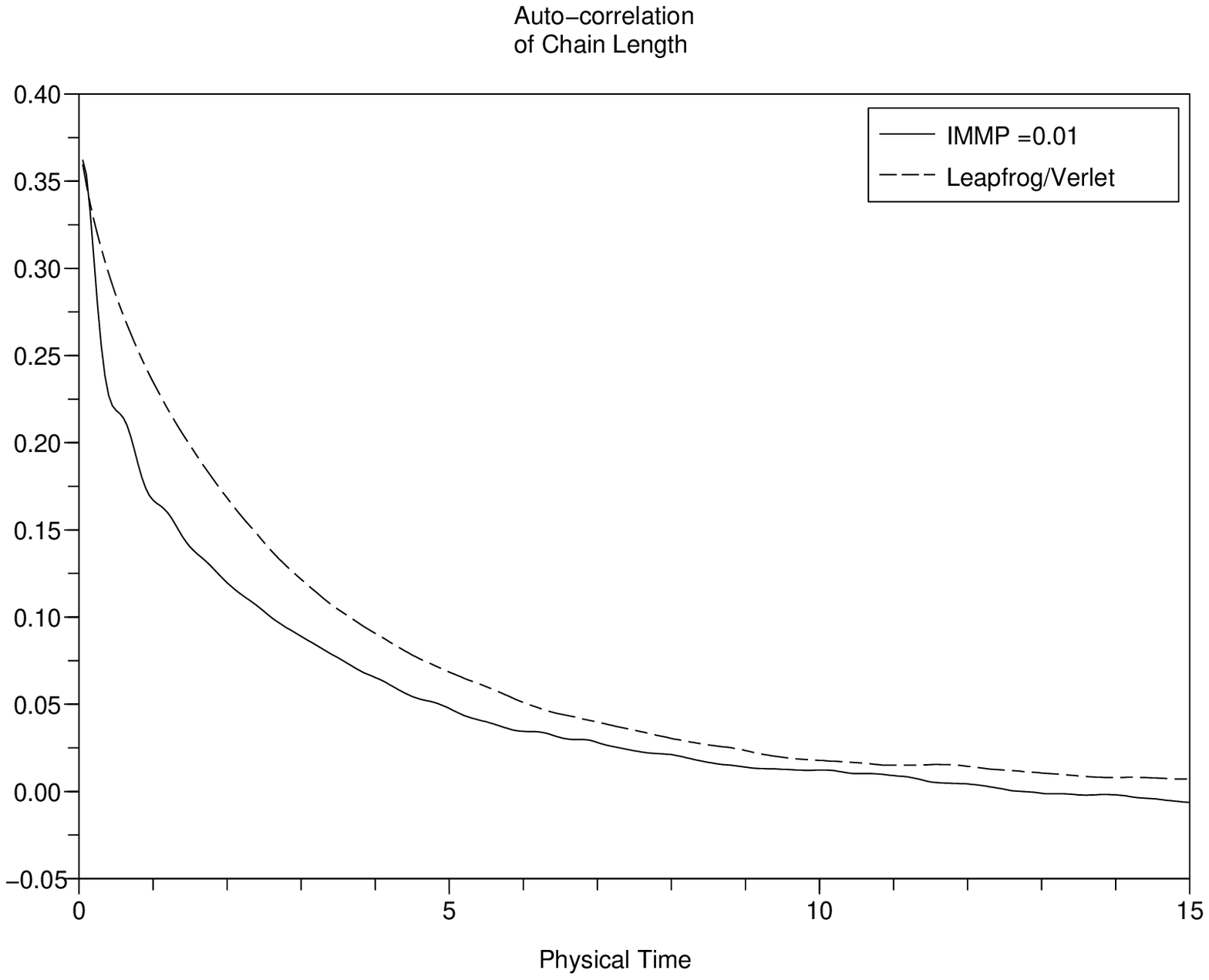}
\includegraphics[width=8.0cm]{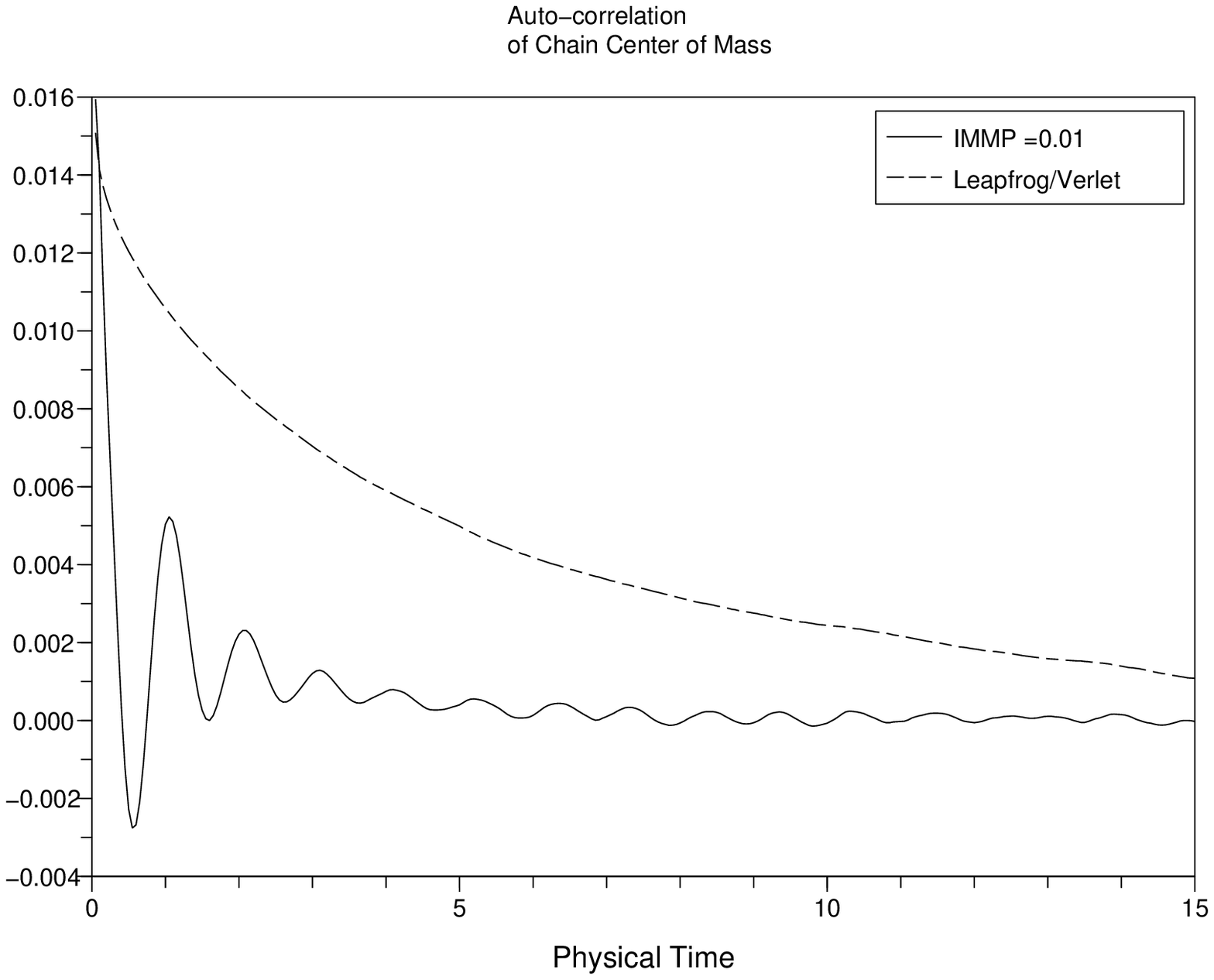}
\vspace{-0.1cm}
\caption{Auto-correlation of a trajectory of the chain length and the center of mass.
        Simulated by the IMMP method  with the penalty $\nusc^2 = 10^{-2}$.}\label{f:corr}
\end{figure}

\medskip
\noindent{\sc Test II:} {\it Relaxation of time-step restrictions.}

\begin{figure}[ht]
  \centerline{%
        \subfigure[Acceptance rate]{\label{f:metro}
         \includegraphics[width=8cm]{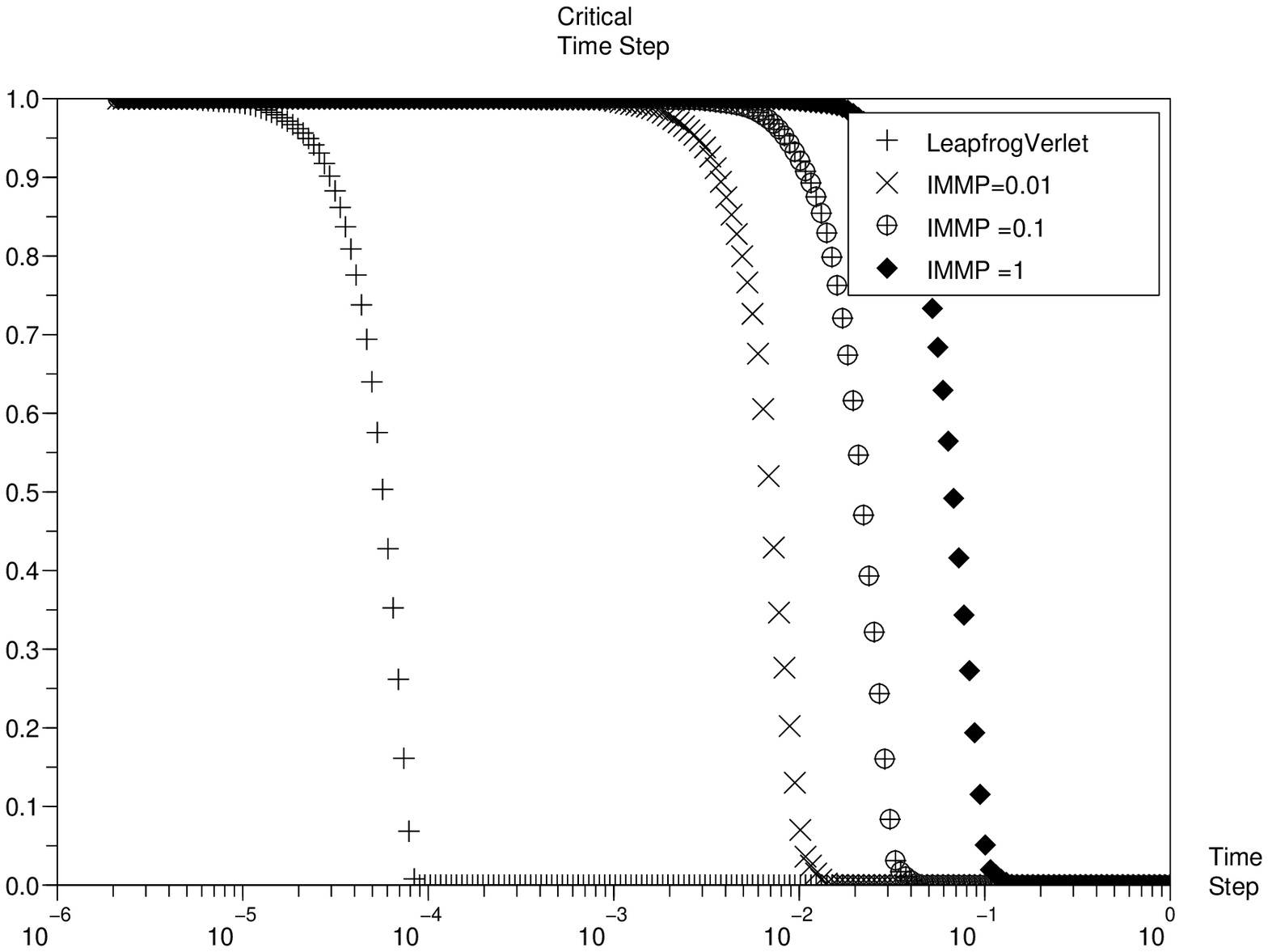}}
        \hspace{0.5cm}
        \subfigure[Critical time-step relaxation]{\label{f:metroscale}
        \includegraphics[width=8cm]{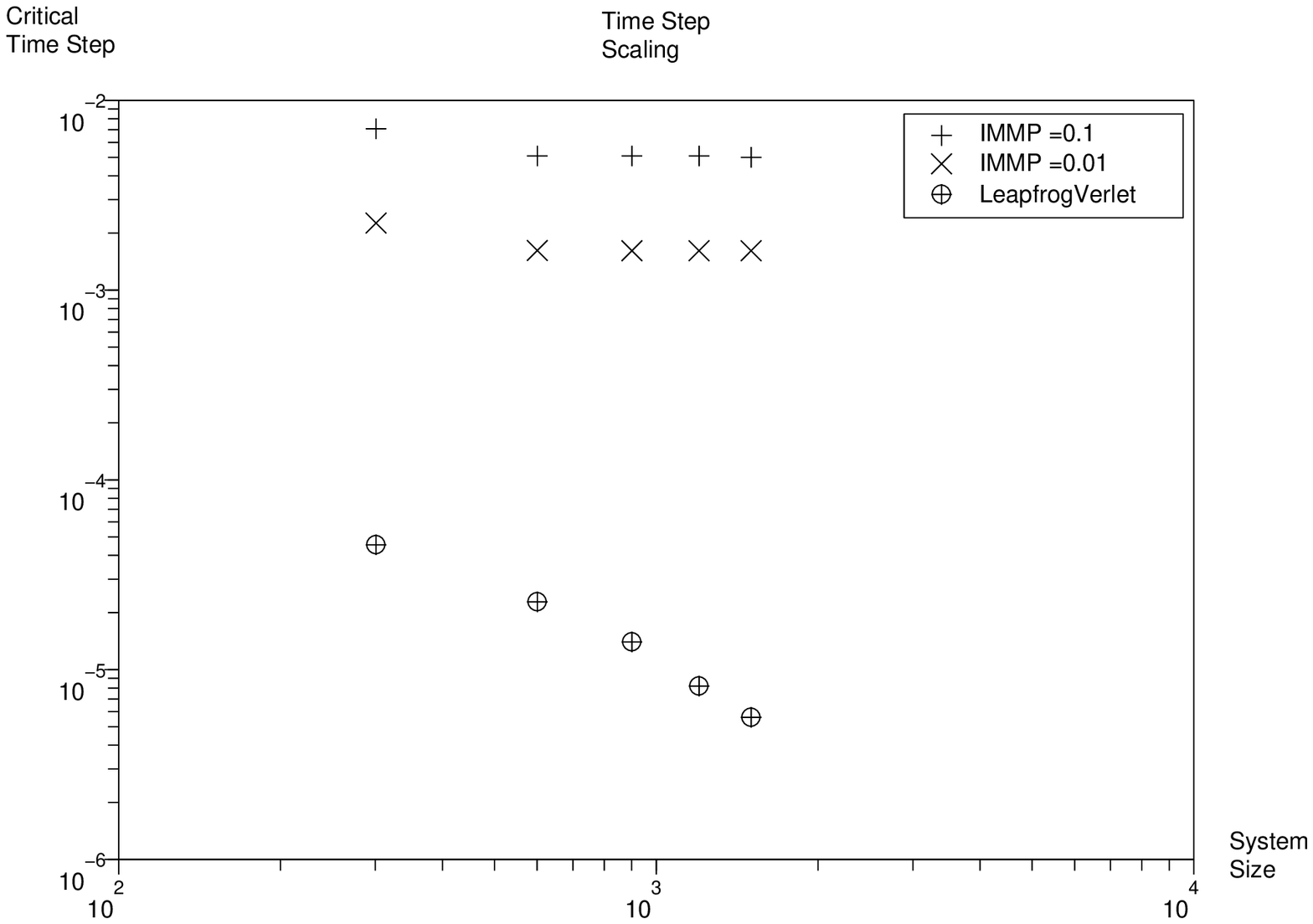}}
    }%
   \caption{(a) The average acceptance rate of the Metropolis step for different values of the
            mass-matrix penalization $\nusc^2= \{1,10^{-1},10^{-2} \}$ and the system size $N=100$.
         (b) Critical time-step that achieves the Metropolis acceptance rate of $0.5$.
             Comparison between the IMMP dynamics ($\nusc^2 = \set{ 10^{-1}, 10^{-2}}$),
             and the exact dynamics, for different values
             of the system size $N$.}\label{f:timestep}
\end{figure}

We first observe in Figure~\ref{f:metro} the relationship between the time-step $\dt$ and the mass penalty $\nusc$.
The average Metropolis acceptance rate is computed for different time-steps  and different values of the
mass-matrix penalization.
The critical time-step, for which the average acceptance rate is approximately $0.5$, is
significantly increased by the mass-matrix penalization, and is broadly proportional to the
re-scaled mass-matrix penalty $\nusc$.

In Figure~\ref{f:metroscale}, the critical time-step is studied for different values
of the chain size $N$ and of the mass-matrix penalization, and it is compared to the exact dynamics. Using linear regression the dependence of the critical time-step on the system size
is estimated and compared to the exact Verlet integration in Table~\ref{t:scale}\footnotemark[1] \footnotetext{Digits inside the brackets indicate the typical error of the last digit.}. These results are in a good agreement with Proposition~\ref{p:stab2} that characterizes
the dependence for the linear case.
\begin{table}[h]
\begin{center}
\begin{tabular}{|c|ccc|}
\hline
  & $\nusc^2 = 10^{-1}$ &  $\nusc^2 = 10^{-2}$  & Leapfrog/Verlet  \\
\hline
 $\dt \propto N^{-\alpha}$ & $\alpha = 0.2 \, (1)$ &$\alpha = 0.2 \,  (1)$ & $\alpha = 1.2  \, (1)$   \\
\hline
\end{tabular}
   \caption{\label{t:scale} Scaling of the critical time-step for the Metropolis acceptance rate $r=0.5$.}
\end{center}
\end{table}

\medskip
\noindent
{\it Conclusions.}

The presented numerical studies demonstrate that for a prescribed macroscopic timescale
of the chain, the IMMP allows for arbitrary large (of order $N$) relaxation of the time-step restriction.
This observation shows that for a large system, $N \to \infty$, the IMMP method can achieve arbitrarily
large computational savings when compared to the standard Verlet algorithm used for integrating the exact dynamics.

\section{Numerical analysis of the particle chain}\label{s:Nanal}
In this section, the particle chain model \eqref{e:HN} is analyzed in the special case of harmonic interactions. 
We consider the thermodynamic limit $N\to +\infty$ where $N$ is the size of the system. It is shown that the macroscopic dynamics of the IMMP method behaves continuously (uniformly with $N$, and in the $L_2$ norm for the position profile) with respect to the re-scaled mass penalty parameter $\nusc$.

At the same time, the time-step stability of the IMMP numerical scheme \eqref{e:Hnum} is compared with the 
standard Verlet scheme, and the critical time step is shown to be increased by a factor $\nusc N$.

From the spectral point of view, the IMMP method behaves in this linear case as a low-pass filter. This proves, in this simplified case, the ability of IMMP method to respect macroscopic dynamical equivalence, while saving computational time up to a factor of order $\BIGO(N)$.

\subsection{Conservation of macroscopic dynamics}

We consider the penalized Hamiltonian  \eqref{e:HNpen} when the interaction potential is quadratic $v_{\INTP}(r) = r^2/2$.
The penalizing matrix is the identity matrix $\MASSTZ = \ID$, hence the penalized 
mass-tensor becomes $\penN{\MASST} = \ID - \nusc^{2}\dl$,
and the fluctuation/dissipation tensor is taken proportional to the identity matrix.
The system of equations \eqref{e:spde} then becomes
\be \label{e:spdelin}
\syst{
& \dot{q} =   ( \ID - \nusc^{2} \dl  )^{-1} p & \\
&\dot{p} =   \dl\,q  - v'_{\EXTP}(q) - \gamma \dot{q} +   \sigma \sqrt{N}\dot{W}  \COMMA&
}
\ee
with fluctuation/dissipation identity $\gamma =  2 \sigma^2\beta^{-1}$. The associated canonical equilibrium distribution is then given by the re-scaled inverse temperature $\beta_N = \beta N^{-1}$.

In order to treat the limit $N\to\infty$, we introduce the $\ell_2$-norm in the position space
\[
\norm{q}_{\ell_2}^2 := \frac{1}{N} \sum_{i=1}^N q_i^2  = \frac{1}{N} q^T q \COMMA
\]
as well as the $h_{-1}$-norm in the momentum space
 \[
 \norm{p}_{h_{-1}}^2 = \norm{(-\dl)^{-1/2}(p-\frac{1}{N}\sum_{i=1}^N p_i)}^2_{\ell_2}+
    \left(\frac{1}{N}\sum_{i=1}^N p_i\right)^2.
 \]
In the above expression, $\tfrac{1}{N}\sum_{i=1}^N p_i$ can be seen as the orthogonal projection in $\ell^2$ on the one dimensional kernel of the Neumann discrete Laplacian $\dl$.  $\norm{q}_0^2 +\norm{p}_{-1}^2$ is a quadratic form that endows the phase-space with a Hilbert space structure.
\bpro[Convergence of the macroscopic dynamics]
Assume that the exterior potential $v_{\EXTP}$ is bounded and that its derivative satisfies the Lipschitz condition
$$
\norm{v'_{\EXTP}(q_2) - v'_{\EXTP}(q_1)}_{h_{-1}} \leq L_v \norm{q_2 - q_1}_{\ell_2} \COMMA
$$
where $L_v$ is independent of $N$. For any $T>0$, let $t \mapsto (p^{\nusc}(t),q^{\nusc}(t))$ be the solution for $t\in [0,T]$ of the evolution equation \VIZ{e:spdelin} with the initial condition
\[
(p^{\nusc}(0),q^{\nusc}(0)) = (\penN{\MASST}^{-1/2} p^0(0),q^0(0)) \COMMA
\]
where $(p^0(0),q^0(0))$ is distributed according to the original equilibrium canonical distribution (associated with \eqref{e:HN}-\eqref{e:beta}). Then for all $t\in [0,T]$ one has the uniform convergence
 \[
 \lim_{\nusc \to 0} \limsup_{N \to + \infty } \EXPECT{\norm{q^{\nusc}(t)-q^{\nusc=0}(t)}^2_{\ell_2}} = 0 \PERIOD
 \]
\epro
\begin{proof} We write $X = (q,p)$,  and introduce the norm:
\[
\norm{X}_{\nusc} = \norm{q}_{\ell_2} + \norm{ \penN{\MASST}^{-1/2} p }_{h_{-1}}\PERIOD
\]
The system \eqref{e:spdelin} becomes a stochastic differential equation in the form
\begin{equation}\label{eqoper}
d X^{\nusc}_t = A_{\nusc} X^{\nusc}_t + F(X^{\nusc}_t) + \Sigma d W_t \COMMA
\end{equation}
where by definition
\[
A_{\nusc} := \bmat
               0                         & (\ID -\nusc^2 \dl )^{-1}  \\
               \dl & 0 \emat \COMMA\;\;
\;\;
F(X) := \bmat
       0 \\-v'_{\EXTP}(q)-  \gamma p
       \emat
        \COMMA\;\;
\Sigma :=  \bmat
       0 \\
       \sqrt{N}\sigma
       \emat\PERIOD
\]
Duhamel formula gives an implicit expression for differences of solutions of \VIZ{eqoper} with the same noise
\begin{eqnarray}\label{e:spdemild}
X^{\nusc}_{t} - X^{0}_{t} &=& \pare{\EXP{A_{\nusc} t} - \EXP{A_{0} t}} X^{0}_0 + \int_{0}^{t}\pare{\EXP{A_{\nusc} (t-s)} - \EXP{A_{0} (t-s)} }\,(F(X^{0}_s)\,ds + \Sigma \,dW_s)  \nonumber \\
&& + \EXP{t A_{\nusc} } (X^{\nusc}_0 - X_{0}^0) +
                  \int_{0}^{t}\EXP{A_{\nusc} (t-s)}( F(X^{\nusc}_{s})-F(X^0_{s})) ds \PERIOD
\end{eqnarray}
We estimate the individual terms on the right hand side  in \VIZ{e:spdemild}.
We define $P$ as the coordinate transformation associated with the orthonormal spectral decomposition
\[
-\dl = P^{-1} \DIAG(\delta_{0},\dots,\delta_{N-1}) P\COMMA
\]
where $P P^{T} = \ID$, and the eigenvalues of the discrete Neumann Laplacian are given, for $k=0,\dots,N-1$, by
\begin{equation}\label{eigenval}
\delta_{k} =
4N^2\sin^2\left(\frac{k\pi}{2N}\right) \limop{\sim}_{N \to \infty} k^2 \pi^2 \PERIOD
\end{equation}
Denoting the spectral coordinates
 \[
\hat{X}=(\hat{q},\hat{p})=(N^{-1/2}Pp, N^{-1/2}Pq)
\]
we have
\[
 \norm{X}^2_{\nusc} = \hat{p}^2_0 + \sum_{k=1}^{N-1} \frac{\delta_k}{1+\nusc^2\delta_k} \hat{p}^2_k +\sum_{k=0}^{N-1} \hat{q}^2_k \PERIOD
\]
The spectral decomposition leads to a block diagonal form of the operator $\EXP{A_{\nusc} t}$ with diagonal
$2\times 2$ blocks in the spectral basis
\[
\EXP{\widehat{A_{\nusc}}^{(0)} t} = \bmat
        1 & t  \\
          0 &  1
         \emat \COMMA
\]
as well as for $k=1,\dots,N-1$
\[
\EXP{\widehat{A_{\nusc}}^{(k)}} = \bmat
         0  &  (1 + \nusc^{2}\delta_k)^{-1}\\
\delta_k         &  0
         \emat \COMMA
\]
where $\widehat{A_{\nusc}}^{(k)}$ is the $2\times 2$ block associated with the coordinates $(\hat{q}_k,\hat{p}_k)$. Since $\widehat{A_{\nusc}}^{(k)}$ conserves the $k$-mode energy $\delta_k  (\hat{q}_k)^2 + (1+\nusc^2 \delta_k)^{-1} (\hat{p}_k)^2$, one can check that for any $N \geq 1$
\[
\normop{\EXP{A_{\nusc} t}}_{\nusc}^2 \leq  2+2t^2 \PERIOD
\]
Similarly, since in the sense of symmetric matrices $\penN{\MASST}^{-1/2} \leq \ID$, we have the bound
\[
\norm{ F(X^{\nusc}) - F(X^0) }_{\nusc} \leq \norm{\penN{\MASST}^{-1/2} \pare{v'_{\EXTP}(q^{\nusc}) -v'_{\EXTP}(q^0) +\gamma (p^{\nusc}-p^0) }}_{h_{-1}} \leq (L_F + \gamma)\norm{X^{\nusc} - X^0 }_{\nusc} \PERIOD
\]
Using random independence of Brownian increments we compute
\begin{eqnarray*}
 \EXPECT{\norm{   \int_{0}^{t} \pare {\EXP{A_{\nusc} (t-s)} -\EXP{A_{0} (t-s)}}\Sigma \,dW_s }^2_{\nusc}} &= &\int_{0}^{t}  \sum_{i=1}^{N} \norm{ (\EXP{A_{\nusc} (t-s)}-\EXP{A_{0} (t-s)}) \Sigma_{.,i}}^2_{\nusc} ds \PERIOD
\end{eqnarray*}
Applying Gronwall lemma in \VIZ{e:spdemild} and collecting all terms we obtain
\begin{equation}\label{e:gronwallN}
\EXPECT{\norm{ X^{\nusc}_{t} - X^{0}_{t}}_{\nusc} ^2} \leq C_T \pare{ \EXPECT{\norm{ X^{\nusc}_{0} - X^{0}_0}_{\nusc} ^2
               + m_T }} \COMMA
\end{equation}
where $C_T$ is independent of $N$, and with $X^0$ being distributed canonically $m_T$ is given
\begin{eqnarray*}
 m_T &=& \sup_{t\in[0,T]} \pare{ \EXPECT{\norm{ \pare{ \EXP{A_{\nusc} t} - \EXP{A_{0} t}} X^0 }_{\nusc}^2}
         + \EXPECT{\norm{ \pare{ \EXP{A_{\nusc} t} - \EXP{A_{0} t}} F(X^0) }_{\nusc}^2}
         + \sum_{i=1}^N  \norm{ \pare{ \EXP{A_{\nusc} t} - \EXP{A_{0} t}} \Sigma_{.,i} }^2_{\nusc} }\PERIOD
\end{eqnarray*}
For a given random vector $X$ such that $\EXPECT{\norm{X}_0^2} < +\infty$, Parseval identity
and the inequality $\norm{\cdot}_{\nusc} \leq \norm{\cdot}_0$ imply
\begin{equation}\label{e:specserie}
\EXPECT{\norm{ \pare{ \EXP{A_{\nusc} t} - \EXP{A_{0} t}} X }_{\nusc}^2} =
  \sum_{k=1}^{N-1} \EXPECT{ \norm{ \pare{\EXP{\widehat{A_{\nusc}} t} - \EXP{\widehat{A_{0}}t} } \hat{X} }_{k,\nusc}^2}
  \leq 2 \sum_{k=1}^{N-1} \EXPECT{\norm{  \hat{X} }_{k,0}^2} \COMMA
\end{equation}
where $\norm{\cdot }_{k,\nusc}$ is the restriction to the $k$-th mode $(\hat{q}_k,\hat{p}_k)$.
Then one has, by orthogonality of $P$,
\[
 \sum_{i=1}^N \EXPECT{\norm{  \widehat{\Sigma_{.,i}}  }_{k,0}^2} =
      \sigma^2  \sum_{i=1}^{N} P^2_{k,i} \frac{1}{\delta_k} \leq \frac{\sigma^2}{\delta_k} \PERIOD
\]
Now the canonical distribution of $X^0=(q^0,p^0)$ can be described as follows. Up to normalization, the distribution of $X^0$ has the density $\EXP{-  \tfrac{\beta}{N} \sum_{i=1}^d v_{\EXTP}(q^0_i)  }$ with respect to the Gaussian distribution with the covariance matrix $\beta^{-1} \ID$ for momenta variables, 
and the covariance matrix $(\beta \dl)^{-1}$ for positions. Thus we have the bound
\[
 \EXPECT{\norm{  \widehat{X^0} }_{k,0}^2} \leq 2\EXP{4\beta \norm{v_{\EXTP}}_{\infty}} \frac{1}{\delta_k \beta} \COMMA
\]
as well as
\[
 \lim_{N \to \infty } \EXPECT{\norm{\widehat{F(X^0)} }_{k,0}^2}
     \leq \EXP{4\beta \norm{v_{\EXTP}}_{\infty}}
     \EXPECT{\norm{ \mathcal{F} \circ F \circ \mathcal{F}^{-1}( \widehat{G^0})}_{k,0}^2} \COMMA
\]
where $\mathcal{F}$ denotes the Fourier series expansion on $[0,1]$ with Neumann conditions, and
$(\widehat{G^0}_k)_{k \geq 1}$ are canonical centered Gaussian i.i.d. variables with the covariance matrix
$\beta^{-1}\bmat 1 & 0\\ 0 & \tfrac{1}{k^2\pi^2} \emat$. By the Lipschitz assumption
the series is bounded
\[
\sum_{k=1}^{+\infty} \EXPECT{\norm{ \mathcal{F} \circ F \circ \mathcal{F}^{-1}( \hat{G^0})}_{k,0}^2}
    \leq (L_v + \gamma) \EXPECT{\norm{\hat{G^0}}_{0}^2} = \sum_{k=1}^{+\infty} \frac{2(L_v + \gamma)}{\beta k^2\pi^2}.
\]
Since $\lim_{\nu \to 0}\normop{\EXP{\widehat{A_{\nusc}}^{(k)} t} - \EXP{\widehat{A_{0}}^{(k)} t}} = 0$,
one can take the limit $N \to +\infty$ and use the uniform convergence of the series in \eqref{e:specserie} to obtain
$\lim_{\nusc \to 0} \lim_{N \to +\infty} m_T =0$ in \eqref{e:gronwallN}. The convergence of the initial condition $\lim_{\nusc \to 0} \lim_{N \to +\infty} \E \norm{ X^{\nusc}_{0} - X^{0}_0}_{\nusc}^2 $ follows by using\ similar arguments. The proof is complete.
\end{proof}

\subsection{Relaxation of time-step stability restriction}
To demonstrate improved stability properties of time integration algorithms we consider the IMMP scheme \eqref{e:immpscheme} associated with the mass-matrix
penalized Hamiltonian \eqref{e:HNpen}. Note that when the constraints are linear, the leapfrog scheme (RATTLE) applied to an implicit
Hamiltonian is identical to the usual leapfrog scheme for the associated explicit Hamiltonian \eqref{e:HNpen}. We restrict the rigorous analysis
to the quadratic interaction potential ($v_{\INTP}(y) = \tfrac{y^2}{2}$),
zero exterior potential ($v_{\EXTP} = 0$), and to the
mass-matrix penalization operator ($\ID-\nusc^2\dl$). The leapfrog scheme is defined as
\begin{equation*}
\left\{\begin{aligned}
p_{n+1/2} &= p_{n} + \frac{\dt}{2} (-\dl) q_{n} \\
q_{n+1} &= q_{n} + \dt \, \penN{M}^{-1} p_{n+1} \\
p_{n+1} &= p_{n+1/2}+ \frac{\dt}{2}  (-\dl)  q_{n+1} \PERIOD
\end{aligned}
\right.
\end{equation*}

Denoting the spectral variables for $k=1,\dots,N-1$
 \begin{equation} \label{e:specvar}
 \left\{\begin{aligned}
\widehat{v}^k & =  \pare{\frac{\delta_k}{1+\nusc^2 \delta_k} }^{1/2} \sqrt{N} P p\\
 \widehat{x}^k & =  \pare{\frac{1+\nusc^2 \delta_k}{\delta_k} }^{1/2} \sqrt{N} P q  \COMMA
 \end{aligned}
 \right.
 \end{equation}
we write
\[
\bmat \widehat{v}_{n+1}^k \\ \widehat{x}_{n+1}^k \emat  = L_k \bmat \widehat{v}_n^k \\ \widehat{x}_n^k \emat \COMMA
\]
where
\[
L_k
=\bmat 1 - \frac{h_k^2}{2} & -h_k + \frac{h_k^3}{4} \\
 h_k & 1 - \frac{h_k^2}{2} \emat
 \COMMA\;\;\;\;\mbox{and}\;\;\;\;h_k = \dt \frac{\delta_k^{1/2}}{\pare{1+\nusc^2 \delta_k}^{1/2}}\PERIOD
\]
Since $\DET\,(L_k)=1$, the standard CFL stability condition is equivalent to \[\abs{\mathrm{Tr}\,(L_k)} \leq 2\]
which is fulfilled  if and only if $h_k \leq 2$ for all $k\leq N-1$.
Thus we arrive at the following bound on the time step
\[
\dt \leq 2 \min_{0\leq k < N} \left(\frac{1+\nusc^2 \delta_k}{\delta_k}\right)^{1/2} \PERIOD
\]

Summarizing the above calculations and recalling \VIZ{eigenval} we have the following characterization
of the stability properties.
\bpro \label{p:stab1} 
Suppose $v_{\EXTP} =0$ and consider a harmonic interaction potential $v_{\INTP}(r) = \tfrac{r^2}{2}$ 
with the mass-matrix penalization $\penN{M}=\ID -\nusc^2 \dl$. 
The leapfrog/Verlet integration of the IMMP harmonic Hamiltonian
\eqref{e:HN} is stable in the spectral sense if and only if
\be \label{e:CFL}
\dt \leq \pare{ 4 \nusc^2 + \frac{1}{ \ds N^2 \\ \sin^2\left(\frac{(N-1)\pi}{2N}\right) } }^{1/2}\PERIOD
\ee
\epro

Since we work with a Metropolis correction of the hybrid Monte-Carlo type, we are also interested in the limiting behavior
of the energy variation compared to the temperature, i.e.,
\[
\beta_N (H(p_{n+1},q_{n+1}) - H(p_{n},q_{n})) \COMMA
 \]
 when
 $(p_{n},q_{n})$
are distributed according to the canonical distribution. This quantity gives the average acceptance rate of
the Metropolis correction. The result we present here is similar to \cite{BesStu}
where the authors analyze infinite dimensional sampling
with the standard Metropolis-Hastings Markov chains.
\bpro \label{p:stab2}
Suppose $v_{\EXTP} =0$ and consider a harmonic interaction potential $v_{\INTP}(r) = \tfrac{r^2}{2}$ with 
the mass-matrix penalization $\penN{M}=\ID -\nusc^2 \dl$.
Suppose the state variable  $X=(\penN{p},q)$ is a random variable distributed according to the canonical distribution
associated with the mass-matrix penalized Hamiltonian \eqref{e:HN}. 
Then the energy variation $\beta_N \Delta H$ after one step
of the leapfrog integration scheme converges in distribution, up to normalization and centering,
to the Gaussian random variable
\[
\frac{\beta_N\Delta H -m_N}{\sigma_N} \xrightarrow[N\to +\infty]{\mathrm{Law}} \mathcal{N}(0,1) \COMMA
\]
with the mean and variance in the infinite size asymptotics for the IMMP method $\nusc > 0 $ and $\dt \equiv\dt_N = o(1)$
\[
  m_N \limop{\sim}_ {N\to +\infty}     \frac{N \dt_N^6}{32 \nusc^{6}} \COMMA \;\;\;\mbox{and}\;\;\;
  \sigma_N^2 \limop{\sim}_ {N\to +\infty} \frac{N\dt_N^6}{16\nusc^{6}} \COMMA
\]
and for the Verlet integration of exact dynamics with $\dt \equiv\dt_N = o(1/N)$
\[
m_N\limop{\sim}_ {N\to +\infty}\frac{5}{8} N^7\dt_N^6  \COMMA\;\;\;\mbox{and}\;\;\;
\sigma_N^2 \limop{\sim}_ {N\to +\infty}  \frac{5}{4}  N^{7} \dt_N^6\PERIOD
\]
\epro
\begin{proof} We start with a canonically distributed state $X=(q,p)$, which is,
by assumption on
the form of the interaction potential,  a Gaussian random vector.
After changing to the spectral coordinates \eqref{e:specvar} we have the spectral representation of the Hamiltonian
\[
\beta_N H = \beta \sum_{k=1}^{N-1}\frac{\delta_k^{1/2}}{2(1+\nusc^2 \delta_k)^{1/2}}\pare{(\widehat{v}^k)^2 + (\widehat{x}^k)^2} \COMMA
\]
and introducing Gaussian random vectors $U$ and $V$ with the identity covariance matrix we can write
\[
\widehat{x}^k = \beta^{-1/2} \frac{(1+\nusc^2 \delta_k)^{1/4}}{\delta_k^{1/4}}U_k\COMMA\;\;\;\mbox{and}\;\;\;
\widehat{v}^k =  \beta^{-1/2} \frac{(1+\nusc^2 \delta_k)^{1/4}}{\delta_k^{1/4}}V_k \PERIOD
\]
We then compute explicitly the change of the Hamiltonian after one step of the leapfrog integration
\begin{equation}\label{deltaH}
\beta_{N} \Delta H = \sum_{k=1}^{N-1} \frac{1}{2} \bmat U_k \\ V_k  \emat^T (L_k^T L_k - \ID ) \bmat U_k \\ V_k  \emat
\PERIOD
\end{equation}
Since $\DET(L_k^T L_k) =1$ the matrix $L_k^T L_k-\ID$ has two positive eigenvalues
$(\lambda_k-1,1/\lambda_k-1)$ which satisfy
\begin{eqnarray*}
&&\lambda_k+1/\lambda_k -2=   \TR(L_k^T L_k-\ID) =\frac{h_k^6}{16} \COMMA \\
&&(\lambda_k-1)^2+(1/\lambda_k -1)^2 =   \TR(L_k^T L_k)^2-2\TR(L_k^T L_k) = \frac{h_k^{12}}{256} + \frac{h_k^6}{8} \PERIOD
\end{eqnarray*}
Combing with \VIZ{deltaH} we find
\[
m_N \equiv \E[\beta_N \Delta H] = \sum_{k=1}^{N-1} \frac{h_k^6}{2^5}\COMMA\;\;\;\mbox{and}\;\;\;
\sigma^2_N \equiv \mathrm{Var}[{\beta_N \Delta H}] = \sum_{k=1}^{N-1}\frac{h_k^6}{2^4}+\frac{h_k^{12}}{2^9}\PERIOD
\]
Moreover, the Lindenberg or simply Lyapunov condition in the general central limit theorem (see \cite{Fel71}) is verified since we work with a sum
of $\chi^2$ random variables, thus concluding the first part of the proof.

Recalling
\[
h_k = \dt \frac{\sin(\frac{k}{N}\frac{\pi}{2})}{(\frac{1}{4N^2}+\nusc^2 \sin^2(\frac{k}{N}\frac{\pi}{2}))^{1/2}} \COMMA
\]
we compute the convergent Riemann sums for $p=6$ and $p=12$.
For the case $\nusc \neq 0$ we have
\[
\lim_{N\to \infty}\frac{1}{N}\sum_{k=1}^{N-1} h^{p}_{k} = \frac{\dt^p}{\nusc^p}\PERIOD
\]
If $\nusc =0$ we obtain
\begin{eqnarray*}
  \lim_{N\to \infty}\frac{1}{N^{p+1}}\sum_{k=1}^{N-1} h^{p}_{k} &=&
  \lim_{N\to \infty}\frac{\dt^p}{N}\sum_{k=1}^{N-1} 2^{p} \sin^{p}\left(\frac{k}{2N}\pi\right) \\
                       &=& \dt^p 2^{p} \int_{0}^{1}\sin^{p}\left(\frac{\pi}{2}x\right)dx\PERIOD
\end{eqnarray*}
Thus for $p=6$ we have that the series sums to $20\dt^6$.
Then the asymptotic behavior follows from the assumption $\dt_N^{12} \ll \dt_{N}^{6}$, and similarly
$\dt_N^{12} N^{12}\ll \dt_{N}^{6} N^6$ in the case $\nusc =0$.
\end{proof}

\brem{\rm
The two propositions proved in this section characterize the restrictions imposed by the
stability of the resulting scheme. 
In Proposition~\ref{p:stab1}, stability in the large system size 
limit, $N \to + \infty$, is equivalent to the inequality \eqref{e:CFL}. In this case the restriction of the
time-step size is imposed by the numerical integrator. On the other hand
the stability for the scheme which uses a Metropolis corrector is linked to the acceptance
rate of the Metropolis step.
In Proposition~\ref{p:stab1}, stability in the large system size limit is equivalent to
the non-vanishing Metropolis acceptance rate, which is equivalent to bounded from above average energy variation 
$m_N$ and bounded variance $\sigma_N$ of  the energy variation.
In either case, the IMMP method ($\nusc>0$) induces a relative increase of order $N$ for the boundary 
of numerical stability as compared to the exact dynamics $\nusc=0$ integrated with the Verlet scheme.
}
\erem

\section{The stiff limit}\label{s:highoscill}
Throughout this section, one introduces a potential function with an explicit dependence with respect to the
fast variables $(q,z) \mapsto U(q,z)$ together with a stiffness parameter
$\epsilon$. The potential energy $V$ can then be written in the form
\begin{equation*}
V(q ) = U(q,\frac{\xi(q)}{\epsilon}) \COMMA
\end{equation*}
with a confining assumption
\[
\inf_{q\in \R^d} U(q,z) \geq K(z) \COMMA
\]
where ${\lim_{z\to \infty} K(z) = +\infty}$ strictly faster than at any logarithmic rate. The functions $\xi$ are then indeed ``fast'' degrees of freedom (fDOFs) in the limit $\epsilon \to 0$, the system being confined to the slow sub-manifold $\SMAN = \set{q \SEP \xi(q)=0}$.

We will prove, that under appropriate scaling of the mass penalty $\nu_{\epsilon} = \tfrac{\nusc}{\epsilon}$, the IMMP method is asymptotically stable in the stiff limit, converging towards standard effective dynamics on the slow manifold $\SMAN$.

\subsection{Thermostatted stiff systems}
The canonical distribution becomes
\be \label{e:boltzmannstiff}
   \mu_{\epsilon}(dp\,dq) = \frac{1}{Z_{\epsilon}}\EXP{ - \beta (\frac{1}{2} p^{T} \MASST^{-1} p +
   U(q,\frac{\xi(q)}{\epsilon})) } dp\, dq \PERIOD
\ee
In the infinite stiffness limit ($\epsilon \to 0$) the measure concentrates on the slow manifold $\SMAN$.
The limit is computed using the co-area formula (see Appendix~\ref{s:measures} for relevant definitions of surface measures).
In order to characterize the limiting measure we introduce the effective potential
\be \label{e:Veff}
U_{\M{eff}}(q) = -\frac{1}{\beta} \ln \int \EXP{- \beta U(q,z) } \, dz \PERIOD
\ee
\blem\label{l:boltzlim}%
In the infinite stiffness limit ($\epsilon \to 0$),
the highly oscillatory canonical distribution \eqref{e:boltzmannstiff} converges $\mu_\epsilon \WEAKLY \mu_0$
(in distribution) towards $\mu_{0}(dp\,dq)$,
 which is supported on $\SMAN$, and defined as
\be \label{e:boltzlim}
\mu_{0}(dp\,dq) = \frac{1}{Z_{0}} \EXP{-\beta (\frac{1}{2} p^{T} \MASST^{-1} p +V_{\M{eff}}(q))} \,dp \, \delta_{\xi(q) = 0}(dq) \PERIOD
\ee
Its marginal distribution in position is given, up to the normalization, by
\be \label{e:boltzlimmarg}
 \EXP{- \beta V_{\M{eff}}(q)}\delta_{\xi(q) = 0}(dq) \PERIOD
 \ee
\elem
\begin{proof}
It is sufficient to consider distributions in the position variable $q$ only. Let $\mathcal{U}^\delta$ be a
$\delta$-neighborhood of $\SMAN$ where
\[
 dq =\epsilon^{\nc} \delta_{\xi(q) = \epsilon z}(dq) \, dz \PERIOD
\]
We construct a decomposition $\ph = \ph_1 + \ph_2$ of continuous bounded observables such that $\supp \ph_1\subset\mathcal{U}^\delta$
and $\supp \ph_2 \cap \mathcal{U}^{\delta/2} = \emptyset$.
Using the confining property of $U(q,\cdot)$ we obtain
\[
\int \ph(q) \EXP{-\beta U(q,\frac{\xi(q)}{\epsilon})} dq =\epsilon^{\nc}
\int \ph_1(q) \EXP{-\beta U(q,z)} \delta_{\xi(q)= \epsilon z}(dq) \, dz  + \BIGO(\EXP{-\beta K(\delta/2\epsilon)})\PERIOD
\]
By continuity of $\epsilon \mapsto \int \ph_1(q) \EXP{-\beta U(q,z)} \delta_{\xi(q)= \epsilon z}(dq)$ and by
the dominated convergence theorem
\[
  \int \ph_1(q) \EXP{-\beta U(q,z)} \delta_{\xi(q)= \epsilon z}(dq) \, dz \to \int \ph_1(q)
  \EXP{- \beta V_{\M{eff}}(q)}\delta_{\xi(q) = 0}(dq)=\int \ph(q) \EXP{- \beta V_{\M{eff}}(q)}\delta_{\xi(q) = 0}(dq)  \COMMA
\]
and the result follows after normalization.
\end{proof}

The infinite stiffness limit ($\epsilon \to 0$) of highly oscillatory dynamics has been studied in a series of papers
\cite{RubUng57, Tak80, Kam85, BorSch97, Rei95, Rei00}. The limiting dynamics can be fully characterized in special cases.
For example, when the highly oscillatory potential is linear and non-resonant
(at least almost everywhere on the trajectory, see \cite{Tak80}), it can be described through adiabatic effective potentials. See also \cite{JahLub06, LebLeg07} for a recent work on some related numerical issues.
However, when the system is thermostatted, one can postulate an ``ad hoc'' effective dynamics (\cite{Rei00}) 
exhibiting the appropriate limiting canonical distribution given by \eqref{e:boltzlimmarg}. 
Such dynamics can be obtained by constraining the system to
the slow manifold $\SMAN$, and adding  a correcting entropic potential, sometimes called ``Fixman'' corrector from \cite{Fix78}, which is due to the geometry of $\SMAN$, and is  given by
\be\label{e:Vfix}
V_{\M{fix}}(q) = \frac{1}{2\beta} \ln \pare{ \mathop{\det} G(q) } \COMMA
\ee
where $G(q)$ is the $\nc \times \nc$ Gram matrix defined in \eqref{e:Gram}.

In general, since the effective potential \eqref{e:Veff} is not explicit, one may need to couple the system with virtual fast degrees
of freedom to enforce the appropriate effective dynamics associated with \eqref{e:Veff}.
The resulting extended Hamiltonian is then defined on the state space $T^{*} \left( \SMAN \times \R^{\nc} \right)$
(the cotangent bundle) and is given by

\begin{equation} \label{e:Heff}
\syst{
&H_{ \M{eff} }(p,p_{z},q,z) =  \frac{1}{2} p^{T} \MASST^{-1} p  +\frac{1}{2} p_{z}^{T} \MASSTZ^{-1} p_{z}
+ U(q,z) + V_{\M{fix}}(q) & \\
&\xi(q) = 0 \PERIOD & \TAG{C}
}
\end{equation}

\begin{Def}[Effective Langevin process with constraints]\label{d:consteff}
The constrained Langevin process associated
with Hamiltonian \eqref{e:Heff} is defined by the following stochastic differential equations
  \begin{equation}\label{e:consteff}
  \syst{
    &\dot{q} =  \MASST^{-1} p &   \\
    &\dot{z} =  \MASSTZ^{-1} p_{z} &  \\
    &\dot{p} = -\Done U(q,z) - \Dq V_{\M{fix}}(q) -\gamma \dot{q} + \sigma \dot{W} - \Dtwo \xi\, \dot{\lambda} &   \\
    &{\dot{p}_{z}} = -\Dz U (q,z)-\gamma_{z} \dot{z} + \sigma_{z} \dot{W}_{z} &\\
    &\xi(q)  = 0\COMMA                 & \TAG{C}
  }\end{equation}
where $\Done$ and $\Dtwo$ are respectively derivatives with respect to the first and second variable of the function $U(q,z)$, $\dot{W}$ (resp. $\dot{W}_{z}$ ) is the standard multi-dimensional white noise, $\gamma$ (resp. $\gamma_{z}$)
  a $d\times d$ (resp. $\nc\times \nc$) symmetric positive semi-definite dissipation matrix,
  $\sigma$ (resp. $\sigma_{z}$) is the fluctuation matrix satisfying $\sigma \sigma^{T} =\tfrac{2}{\beta} \gamma$
  (resp. $\sigma_{z} \sigma_{z}^{T} =\tfrac{2}{\beta} \gamma_{z}$).
  The processes $\lambda \in \mathbb{R}^{\nc}$ are Lagrange multipliers associated with the constraints
  $(C)$ and adapted with respect to the white noise.
\end{Def}

We formulate reversibility of this process as a separate lemma.
\blem
The process defined in \VIZ{e:consteff} is reversible with respect to the associated canonical distribution whose marginal distribution in $(q,p)$ variables is
\begin{equation}\label{e:mueff}
\mu_{\M{eff}}(dp \, dq) = \frac{1}{Z_{\M{eff}}}\EXP{-\beta \pare{ \frac{1}{2} p^{T} \MASST^{-1} p + V_{\M{eff}}(q) + V_{\M{fix}} (q)} } \sigma_{T^{*}\SMAN}(dp\,dq)
\end{equation}
with the $q$-marginal
\[
 \EXP{- \beta V_{\M{eff}}(q)}\,\delta_{\xi(q) = 0}(dq) \PERIOD
\]
When $\gamma$ and $\gamma_{z}$ are strictly positive definite, the process is ergodic.
\elem
\begin{proof}
The process \VIZ{e:consteff} is a Langevin process with mechanical constraints, exhibiting reversibility properties
with respect to the associated Boltzmann canonical measure (see the summary in Appendix~\ref{s:langevinapp}).
Then the $q$-marginal is obtained by remarking that the integration of any function of
$\tfrac{1}{2} p^{T} \MASST^{-1} p + \tfrac{1}{2} p_{z}^{T} \MASSTZ^{-1} p_{z}$ with respect to $d p_z \, \sigma_{T^{*}_q\SMAN}(dp)$
results in a constant independent of $q$.
\end{proof}

The properties of thermostatted highly oscillatory systems are summarized in Table~\ref{t:sum}.
\begin{table}[H]
\begin{center}
\begin{tabular}{|c|cccc|}
\hline
&   \multicolumn{2}{c|}{Finite stiffness} &      \multicolumn{1}{c|}{ Infinite stiffness limit }    &   \multicolumn{1}{c|}{ Infinite stiffness  }            \\
&   \multicolumn{2}{c|}{ $\epsilon >0$  }  &      \multicolumn{1}{c|}{$\epsilon \to 0$}  & \multicolumn{1}{c|}{$\epsilon = 0$}                       \\
\hline
Dynamics & Highly & & \multicolumn{1}{c|}{Adiabatic }      & Effective with \\
        &  oscillatory  &                             & \multicolumn{1}{c|}{(if non-resonant)}    &     constraints \\
         &  $+$ fluct./diss.   &                               &\multicolumn{1}{c|}{ $+$ non-Markov fluct./diss.}   &   $+$ fluct./diss.            \\

\hline
Statistics & Canonical  &  &  \multicolumn{1}{c|}{Positions on $\mathcal{M}_{0}$,} & Canonical on $T^{*}\mathcal{M}_{0}$,  \\
           &            &                                &  \multicolumn{1}{c|}{free velocities.}     &  geometric corrector. \\
\hline
Numerics & Leapfrog/Verlet  &      &    \multicolumn{1}{c|}{Time-step  }                   &      Leapfrog/Verlet with     \\
     &     $+$ fluct./diss.       &                 &    \multicolumn{1}{c|}{restrictions ($\dt = o(\epsilon)$)}               &       constraints     \\
      &         &                      &    \multicolumn{1}{c|}{}               &        $+$ fluct./diss.           \\
\hline

\end{tabular}
\caption{\label{t:sum} Stiff  Hamiltonian systems and associated commonly used numerical methods
($\SMAN$ denotes the slow manifold).
Two different schemes are required
for the stiff system and its effective Markovian approximation.}
\end{center}
\end{table}

\subsection{Stability of the IMMP dynamics}\label{s:infstiff}
We assume that the mass-matrix penalty parameter $\nu\equiv\nueps$ grows to infinity in such a way that
\[
\lim_{\epsilon\to 0} \epsilon \nu_\epsilon = \nusc \PERIOD
\]

The original Hamiltonian with the stiffness parameter is expressed explicitly as
\be\label{e:Heps}
H_{\epsilon}(p,q) =
\frac{1}{2} p^{T} \MASST^{-1} p    + U(q,\frac{\xi(q)}{\epsilon}) \COMMA
\ee
and including the mass-matrix penalization one gets
\be\label{e:Hnueps}
	\peneps{H}(\peneps{p},q) = \frac{1}{2} \peneps{p}^{T} \peneps{\MASST}^{-1} \peneps{p}
                   + U(q,\frac{\xi(q)}{\epsilon}) + V_{\M{fix},\nueps}(q)  \COMMA
\ee
or in its implicit formulation
\begin{equation}\label{e:Himmpeps}\syst{
   &H_{\M{IMMP}}(q,z,p,p_{z}) = \frac{1}{2} p^{T} \MASST^{-1} p  + \frac{1}{2} p_{z}^{T} \MASSTZ^{-1} p_{z} +
   U(q,\frac{z}{\nueps\epsilon} ) + V_{\M{fix},\nueps}(q)\COMMA & \\
   &\xi(q) = \frac{1}{\nueps} z\PERIOD & \TAG{C_{\nueps}}
}\end{equation}
One immediately sees that $H_{\M{IMMP}}$ is non-singular when $\epsilon\to 0$ and converges to the effective Hamiltonian on the slow manifold,
\begin{equation}\label{e:Himmpstab}
\syst{
     &  H_{\M{eff},\nusc}(q,z,p,p_{z}) = \frac{1}{2} p^{T} \MASST^{-1} p  + \frac{1}{2} p_{z}^{T} \MASSTZ^{-1} p_{z}
                                          + U(q,\frac{z}{\nusc}) + V_{\M{fix}}(q) &\\
       &\xi(q) = 0 \PERIOD &\TAG{C} \PERIOD
}
\end{equation}
The expression \eqref{e:Himmpeps} represents a minor generalization of  $H_{\M{eff}}$ \eqref{e:Heff},
but it leads to the same canonical marginal distribution $\mu_{\M{eff}} (dp \, dq)$ in $(p,q)$ variables as given by \eqref{e:mueff}.
The continuity in $\epsilon$ of $H_{\M{IMMP}}$ implies stability of the associated dynamics and their numerical integrators.
We  first derive the limits of the original and penalized canonical distribution.
\bpro[Limits of canonical distributions]
Consider the canonical distributions $\peneps{\mu}(dp \, dq$ associated with the mass penalized Hamiltonian \eqref{e:Hnueps}, but considered with respect to the $(p=\MASST\peneps{\MASST}^{-1} \peneps{p},q)$ variables. In the sense of weak convergence of measures we have $\peneps{\mu}\WEAKLY \mu_{\mathrm{eff}}$ as $\epsilon\to 0$ with
$\mu_{\M{eff}}$ defined by \eqref{e:mueff}.
\epro
\begin{proof}
The first convergence is  proved in  Lemma~\ref{l:boltzlim}. For the second one,  the following notation will be used
\[
\syst{
&\delta_{q,\epsilon z}(dq)  = \delta_{\xi(q) = \epsilon z}(dq)  &\\
&\delta_{p,\epsilon p_z}(dp)  =\delta_{p^T \MASST^{-1} \nabla \xi(q)  = \epsilon\MASSTZ^{-1} p_z}(dp) \PERIOD&
}
\]
To prove the convergence towards $\mu_{\mathrm{eff}}$
we consider a $\delta$-neighborhood $\mathcal{U}^\delta$ of $\SMAN$ where
\begin{eqnarray*}
  dp \, dq &=&\frac{\epsilon^{2\nc}}{\DET\MASSTZ} \delta_{q,\epsilon z}(dq) \, dz \;
  \delta_{p,\epsilon p_z}(dp) \, dp_z \COMMA \\
\end{eqnarray*}
and a decomposition of the bounded observable (in $(p,q)$ variables) $\ph = \ph_1 + \ph_2$
such that $\supp \ph_1\subset\mathcal{U}^\delta$ and $\supp\ph_2\cap\mathcal{U}^{\delta/2} = \emptyset$.
Thus, keeping in mind that $\peneps{p} = \peneps{\MASST} \MASST^{-1} p$, and using the confining property of the potential $U(q,\cdot)$ we obtain
\begin{equation}\label{testfun}
\int \ph(p,q) \EXP{-\beta \peneps{H}} d \peneps{p} \, dq = \int \ph_1(p,q) \EXP{-\beta \peneps{H}} d \peneps{p} \, dq +
 \BIGO(\EXP{-\beta K(\delta/\epsilon)}) \equiv I_\epsilon +  \BIGO(\EXP{-\beta K(\delta/2\epsilon)}) \PERIOD
\end{equation}
Applying the change of variables $\peneps{p} = \peneps{\MASST} \MASST^{-1} p $ yields
\[
d\peneps{p} = \DET(\peneps{\MASST }\MASST^{-1})\, dp
= \nueps^{2\nc} \,\DET{\MASSTZ} \,\DET(G+\frac{1}{\nueps^2}\MASSTZ^{-1})\, dp \COMMA
\]
and setting $  \epsilon\MASSTZ^{-1} p_z = \nabla_q \xi\,\MASST^{-1} p $ and
$\epsilon z = \xi(q) $ we get
\[
\peneps{H}(\peneps{p},q) = \frac{1}{2} p^{T} \MASST^{-1} p    + \nueps^2 \epsilon^2
p_z^{T} \MASSTZ^{-1} p_z + U(q,\frac{z}{\nueps\epsilon}) + V_{\M{fix},\nueps}(q)
= H_{\M{IMMP}}(q,z,p,p_z)\PERIOD
\]
Thus substituting back to \VIZ{testfun} we obtain
\[
 I_\epsilon = (\nueps \epsilon)^{2 \nc} \int \ph_1 \EXP{-\beta H_{\M{IMMP}}(q,z,p,p_z)} \DET(G+\frac{1}{\nueps^2}\MASSTZ^{-1}) \, \delta_{p,
\epsilon p_z}(dp) \, dp_z\, \delta_{q,\epsilon z}(dq) \, dz \COMMA
\]
and thus
\[
I_\epsilon \xrightarrow[\epsilon \to 0]{}  \nusc^{2\nc}
\int \ph_1 \EXP{-\beta H_{\M{eff},\nusc}(q,z,p,p_z)}  \DET(G) \, \delta_{p,\epsilon p_z}(dp) \, dp_z \,
\delta_{\xi(q) = 0}(dq) \, dz
\PERIOD
\]
Using the co-area formula we obtain
\[
\DET(G)\, \delta_{\nabla \xi(q)\MASST^{-1} p  = 0}(dp)\delta_{\xi(q) = 0}(dq) = \sigma_{T^* \SMAN}(dp \, dq) \COMMA
\]
which leads to the final result after integration of the $(p_z,z)$ variables and normalization.
\end{proof}
\brem {\rm
Due to the fast oscillations, the distribution of impulses in the limiting distribution $\mu_{0}$ in
\eqref{e:boltzlim} is uncorrelated, whereas after the mass-matrix penalization, the limiting distribution \eqref{e:penboltz}
has almost surely co-tangent impulses (i.e., satisfying the constraints $\nabla_q \xi \MASST^{-1} p =0$).
This explains the role of the corrected potential energy $V_{\M{fix}}$ taking into account the curvature of $\mathcal{M}_{0}$.
}
\erem

In the next step we inspect the infinite stiffness asymptotic of the penalized dynamics.
\bpro[Infinite stiffness limit]
	When $\epsilon \to 0$ with $\nu\equiv\nueps \sim \tfrac{\nusc}{\epsilon}$ and
        $V(q,\xi(q))= U(q,\tfrac{\xi(q)}{\epsilon})$,
        the IMMP Langevin stochastic process \eqref{e:IMMP} converges weakly towards
        the following coupled limiting processes with constraints
	\begin{equation}\label{e:limdyn}
        \syst{
    	  &\dot{q} =  \MASST^{-1} p \COMMA &  \\
    	  &\dot{p} = -\Done U(q,\frac{z}{\nusc}) - \Dq V_{\M{fix}}(q) -\gamma \dot{q} + \sigma \dot{W} - \Dq\xi \dot{\lambda}\COMMA  & \\
    	  &\xi(q)  =  0  \COMMA                                        &  \TAGG{C} \\
    	  &\dot{z} =  \MASSTZ^{-1} p_{z} \COMMA &  \\
    	  &{\dot{p}_{z}} = -\frac{1}{\nusc}\Dtwo U(q,\frac{z}{\nusc})-\gamma_{z} \dot{z} + \sigma_{z} \dot{W}_{z} \PERIOD &
	}
        \end{equation}
	where $\Done$ and $\Dtwo$ are respectively derivatives with respect to the first and second variable of the function $U(q,z)$,
        and $\{\lambda_t\}_{t\geq 0}$ are adapted stochastic processes defining the Lagrange multipliers associated with the constraints $(C)$.

	The process $\{q_t,p_t\}_{t\geq 0}$ defines an effective dynamics with constraints
        (Definition~\ref{d:consteff}) for thermostatted highly oscillatory systems.
        It is reversible with respect to its stationary canonical distribution given by $\mu_{\mathrm{eff}}$ \eqref{e:mueff},
        and is ergodic when $(\gamma,\gamma_z)$ are strictly positive definite.
\epro
\begin{proof} The proof is similar to the proof of Proposition~\ref{p:largepen}. Here we have
\[
 \Dq U = \Done U + \frac{1}{\epsilon} \Dq \xi \Dtwo^T U \COMMA
\]
and \eqref{e:IMMP} translates, up to a change of Lagrange multipliers, into
\begin{equation}\syst{
    &\dot{q} =  \MASST^{-1} p  &  \\
    &\dot{z} =  \MASSTZ^{-1} p_{z} &  \\
    &\dot{p} = -\Done U - \Dq V_{\M{fix},\nu}(q) -\gamma \dot{q} + \sigma \dot{W} - \Dq \xi \,\dot{\lambda} & \\
    &{\dot{p}_{z}} = - \frac{1}{\nueps \epsilon} \Dtwo U-\gamma_{z} \dot{z} + \sigma_{z} \dot{W}_{z} +  \frac{\dot{\lambda}}{\nueps} & \\
    &\xi(q) = \frac{z}{\nueps}\COMMA & \TAG{\peneps{C}} \PERIOD
  }\end{equation}
The rest follows the proof of Proposition~\ref{p:largepen}.

\end{proof}
\brem {\rm
When $\nusc \to +\infty$, by a classical averaging argument (see, e.g., \cite{Kam85}), one can check that the limiting dynamics
are the effective dynamics pointed out in \cite{Rei00}
\begin{equation} \label{e:avlimdyn}\syst{
    	&\dot{q} =  \MASST^{-1} p & \\
    	&\dot{p} = -\Dq U_{\M{eff} } (q) - \Dq V_{\M{fix}}  -\gamma \dot{q} + \sigma \dot{W} - \Dq\xi\, \dot{\lambda}& \\
    	&\xi(q)  = 0  \PERIOD & \TAG{C}
}\end{equation}
with the stationary canonical distribution \VIZ{e:mueff}.
}
\erem

\subsection{Stability of the IMMP integrator}
The numerical scheme (Scheme~\ref{d:scheme} proposed for the IMMP method \eqref{e:IMMP}) is also stable in the limit of infinite stiffness
$\epsilon \to 0$.
Recall that we consider a reversible, measure preserving numerical flow
$\Phi_{\dt}^{\nueps}(p,p_z,q,z)$ associated with Hamiltonian \eqref{e:Himmpeps}  $H_{\M{IMMP}}$ with constraints (modified potentials could similarly be considered).
\bpro[Asymptotic stability]
In the limit $\epsilon\nueps \to \nusc$, the numerical flow $\Phi_{\dt}^{\nueps}$ associated with the leapfrog/Verlet integrator with constraints
for the IMMP Hamiltonian \eqref{e:Himmpeps} converges towards the numerical flow
$\Phi_{\dt}^{\nusc}$, which is the leapfrog/Verlet integrator with geometric constraints associated with effective Hamiltonian \eqref{e:Himmpstab} on the slow manifold.
\epro
\begin{proof}
The statement is a direct consequence of the implicit function theorem and  the continuity of the leapfrog integrator
with constraints \eqref{e:immpscheme} with respect to the parameter $\nusc=\epsilon \nueps$.
Indeed, considering the shift of Lagrange multipliers
\[
\lambda \to \lambda +\frac{1}{\epsilon} \Dtwo U
\]
and taking the limit $\epsilon \to 0$ we obtain the appropriate leapfrog scheme
\begin{equation*}\syst{
		&p_{n+1/2} = p_{n} -  \frac{\dt}{2} \Done U(q_{n},\frac{z_n}{\nusc}) - \Dq \xi(q_{n})  \lambda_{n+1/2} &\\
		&p_{n+1/2}^z = p_{n}^z -  \frac{\dt}{2 \nusc} \Dtwo U(q_{n},\frac{z_n}{\nusc})     & \\
		&q_{n+1} =  q_{n}+ \dt \MASST^{-1} p_{n+1/2}            & \\
		&z_{n+1} = z_n + \dt \MASSTZ^{-1} p_{n+1/2}^z             & \\
		&\xi(q_{n+1}) = 0                       & \TAG{C_{1/2}} \\
		&p_{n+1} = p_{n+1/2} -  \frac{\dt}{2} \Done U(q_{n+1},\frac{z_{n+1}}{\nusc}) - \Dq \xi(q_{n+1})  \lambda_{n+1} &\\
		&p_{n+1}^z = p_{n+1}^z -  \frac{\dt}{2 \nusc} \Dtwo U(q_{n+1},\frac{z_{n+1}}{\nusc})  & \\
		&\Dq \xi(q_{n+1})\MASST^{-1}p_{n+1} = 0 \PERIOD & \TAG{C_1}
}\end{equation*}
\end{proof}

By convergence of the Hamiltonian \eqref{e:Himmpeps} to \eqref{e:Himmpstab},  similar asymptotic stability properties holds when a Metropolis step is introduced.

The results and properties discussed in this section are summarized in Table~\ref{t:sum2}.
\begin{table}[H]
\begin{center}
\begin{tabular}{|c|cccc|}
\hline
& Zero mass    &  \multicolumn{2}{|c|}{Positive }          &    \multicolumn{1}{|c|}{ Infinite }  \\
& penalization &  \multicolumn{2}{|c|}{ mass-penalization} &    \multicolumn{1}{|c|}{ stiffness limit }  \\

&$\nu = 0$              &  \multicolumn{2}{|c|}{ $\epsilon,\nu >0$  } &    \multicolumn{1}{|c|}{ $\epsilon\to 0$, $\tfrac{\nu}{\epsilon} \to \nusc$}               \\
\hline
Dynamics & Highly oscillatory & IMMP && Effective with  \\
         & $+$ fluct./diss.   &  $+$ fluct./diss.     &                                &    constraints$+$ fluct./diss.    \\

\hline
Statistics  &Canonical&   Canonical with      &                &  Canonical on \\
            &         & correlated velocities &                &      $T^* \SMAN$                 \\
\hline
Numerics & \multicolumn{4}{c|}{ IMMP $+$ fluct./diss. }     \\
\hline
\end{tabular}

\caption{\label{t:sum2} The IMMP dynamics and the Verlet numerical integration are both asymptotically stable in the infinite stiffness
regime if $ \tfrac{\nu}{\epsilon} \to \nusc < +\infty$.
If the mass-penalization vanishes ($\nu =0$) one recovers the original physical stiff system.
The canonical distribution is always exact in the position variable. Notice that due to the penalized mass-matrix ($\nu>0$)
the statistics have correlated velocities.}
\end{center}
\end{table}

\appendix

\section{Surface measures}\label{s:measures}
Let $\mathbb{R}^{d}$ be endowed with the scalar product given by the positive definite matrix $\MASST$, and consider
$\SMANZ{z}$ a family
of sub-manifolds of co-dimension $\nc$ implicitly defined by $\nc$ independent functions
$\SMANZ{z} = \{ x \in \mathcal{M} \SEP \xi_{1}(q)=z_{1},..,\xi_{\nc}(q)=z_{\nc} \}$ for $z$ in a neighborhood of the origin.
For each $z$ in a neighborhood of the origin the \emph{conditional measure} $\delta_{\xi(q)=z}(dq)$ is a measure on $\SMANZ{z}$ defined
in such a way that it satisfies the chain rule for conditional expectations with respect to the Lebesgue measure $dq$,
i.e.,
\be
dq = \delta_{\xi(q) = z}(dq) \, dz\PERIOD
\ee
The {\em surface measure} $\sigma_{T^{*}_{q}\SMANZ{z}}(dp)$  is the Hausdorff measure induced by the metric $\MASST^{-1}$
on the co-tangent space $T^{*}_{q}\SMANZ{z}=\set{ p \SEP \nabla_q^T \xi(q)\MASST^{-1} p  = 0}$; and in the same way,
$\sigma_{\SMANZ{z}}(dq)$  is the  Hausdorff measure induced by the metric
$\MASST$ on the sub-manifold $\SMANZ{z}$. It is important to note that, although this is not explicit in notation, $\sigma$
is defined with respect to the mass-tensor $\MASST$ of the mechanical system. The Liouville measure $\sigma_{T^{*}\SMANZ{z}}(dp \, dq)$
on the co-tangent bundle
$T^{*}\SMANZ{z}$ is the volume form induced on
\[
T^{*}\SMANZ{z} = \set{(p,q)\SEP \nabla_q^T \xi(q)\MASST^{-1} p  = 0, \quad\xi(q)=z }
\]
by the usual symplectic form $dp \wedge dq$. It can be described in terms of surface measures as follows
\[
\sigma_{T^{*}\SMANZ{z}}(dp \, dq) = \sigma_{T^{*}_q\SMANZ{z}}(dp) \, \sigma_{\SMANZ{z}}(dq) \PERIOD
\]

Finally, the co-area formula (see \cite{Fed69} for a general reference) defines the relative probability density between
$\delta_{\xi(q)=z}(dq)$ and $\sigma_{\SMANZ{z}}(dq)$.
\begin{Pro}[Co-area formula]
Given the invertible Gram matrix associated with the constraints $\xi(q) =z$ in a neighborhood of $\SMANZ{z} = \{ q \SEP \xi(q)= z \}$
\[
G(q) = \nabla_q^T \xi\, \MASST^{-1} \nabla_q \xi\COMMA
\]
one has
\[
\delta_{\xi(q)=z}(dq) =  \frac{1}{\sqrt{ \DET G(q)}} \sigma_{\SMANZ{z}}(dq) \PERIOD
\]
\end{Pro}
\section{Langevin processes}\label{s:langevinapp}
Defining the Poisson bracket
\[
\poisson{\ph_1,\ph_2} = \Dp^T\ph_1  \Dq \ph_2  - \Dp^T\ph_2  \Dq \ph_1 \COMMA
\]
and the dissipation tensor
\[
\diss(q) = \sigma q \COMMA
\]
where $\sigma$ is the fluctuation matrix in Definition~\ref{d:langevin},
the Markov generator of the Langevin process in Definition~\ref{d:langevin} is
\[
\mathcal{L} = \poisson{\,\cdot\, , H} + \frac{1}{\beta}
\poisson{\diss,\poisson{\diss^T, \, \cdot \, }\EXP{- \beta H}}\EXP{\beta H} \PERIOD
\]
The generator $\mathcal{L}$ satisfies
\[
\int \ph_1 \,\mathcal{L} (\ph_2) \EXP{- \beta H} \,dp \, dq =
\int \mathcal{L}^* (\ph_1) \,\ph_2 \,\EXP{- \beta H} \,dp \, dq \COMMA
\]
where
\[
\mathcal{L}^* = \poisson{\,\cdot\, , -H} + \frac{1}{\beta}
\poisson{\diss,\poisson{\diss^T, \, \cdot \, }\EXP{- \beta H}}\EXP{\beta H} \PERIOD
\]
The generator $\mathcal{L}^*$ defines a Langevin process with the time-reversed Hamiltonian ($-H$).
Reversibility of the process implies that the canonical measure is stationary.
Furthermore, if the initial state of the system is a canonically distributed random variable, the probability distribution
of a trajectory after the time-reversal is given by a Langevin process with the generator $\mathcal{L}^*$.
When $H$ has the form $H(p,q) =  \tfrac{1}{2} p^{T} \MASST^{-1} p    + V(q) $, reversal of impulses ($p\to -p$) leads
to time-reversed dynamics, and a process with
generator $\mathcal{L}^*$ can be constructed by the following simple steps:
 \begin{enumerate}
 \item Reverse momenta ($p\to -p$).
 \item Draw a random path with generator $\mathcal{L}$.
 \item Reverse again momenta ($p\to -p$).
 \end{enumerate}
When holonomic constraints, for instance, of the form
\[
\Xi(p,q) = \zeta \Leftrightarrow \syst{
p^{T} M^{-1} \nabla_{q} \xi &=& 0 \\
\xi(q) &=& z
}
\]
are introduced, it is useful to define the Poisson bracket on the co-tangent bundle $T^*\SMANZ{z}$
\[
\poisson{\ph_1,\ph_2}_{\SMANZ{z}} = \poisson{\ph_1,\ph_2} - {\ds \sum_{a,b}}\poisson{\ph_1,\Xi^a} \Gamma^{-1}_{a,b} \poisson{\Xi^b,\ph_2} \COMMA
\]
where $\Gamma$ is the symplectic Gram matrix of the full constraints
\[
\Gamma^{a,b} = \poisson{\Xi^a,\Xi^b} \PERIOD
\]
As a basic result of symplectic geometry (see \cite{Arn89}), one recovers the divergence formula with respect to the bracket
$\poisson{\,\cdot \,\,, \cdot\,}_{\SMANZ{z}} $ and the Liouville measure $\sigma_{T^{*}\SMANZ{z}}(dp \, dq)$
\[
\int \poisson{\,\cdot \,, \cdot\,}_{\SMANZ{z}} \sigma_{T^{*}\SMANZ{z}}(dp \, dq) = 0 \PERIOD
\]
Given a constrained Langevin process in a stochastic differential equation form
\begin{eqnarray*}
    \dot{q} &=& \Dp H \COMMA \\
    \dot{p} &=& -\Dq H -\gamma \dot{q} + \sigma \dot{W} - \Dq \xi \,\dot{\lambda}\COMMA
   \end{eqnarray*}
where $\lambda$ are Lagrange multipliers associated with the constraints $\xi(q)=0$,
adapted with respect to the noise $\dot{W}$, the process $\{p_t,q_t\}_{t\geq 0}$ obeys hidden 
velocity constraints and is characterized
by the stochastic differential equations
\begin{eqnarray*}
    \dot{q} &=& \Dp H  +  \Dp \Xi \,\dot{\Lambda}\COMMA \\
    \dot{p} &=& -\Dq H -\gamma \dot{q} + \sigma \dot{W} - \Dq \Xi\,\dot{\Lambda}\COMMA
\end{eqnarray*}
where $\Lambda$ are Lagrange multipliers associated with the full constraints $\Xi(p,q)=0$.
The Markov generator of this process can be written in the form
\[
\mathcal{L}_{\SMANZ{z}} = \poisson{\,\cdot \,, H}_{\SMANZ{z}} + \frac{1}{\beta}
\poisson{\diss,\poisson{\diss^T, \, \cdot \, }_{\SMANZ{z}}
\EXP{- \beta H}}_{\SMANZ{z}}\EXP{\beta H} \COMMA
\]
demonstrating the reversibility with respect to the constrained canonical measure
$\EXP{-\beta H} \sigma_{T^{*}\SMANZ{z}}(dp \, dq)$.

\section{Exact sampling of fluctuation/dissipation perturbations}\label{s:exactflucdiss}
In this section, we recall how to perform exact sampling of fluctuation/dissipation perturbations.
Since we only work with impulses, we refer to the system by using the impulse variables $p$ only.
Note that throughout the paper, we also use extended variables $(p,p_z)$,
however, the presentation that follows covers general cases.
The kinetic energy of the system is $\tfrac{1}{2}  p^{T}\MASST  p$. We impose constraints
$p^{T} \MASST^{-1} \Dq \xi = 0$ on impulses, thus $p\in T^{*}_{q}\mathcal{M}$
and hence the associated orthogonal projector on $T^{*}_{q}\mathcal{M}$ is
\[
P  = \ID - \Dq \xi \,G^{-1} \Dq^T \xi \,\MASST^{-1} \PERIOD
\]
The stochastic differential equations of motion on impulses that are integrated on a time-step interval are
\begin{equation}\label{e:sdeimp}\syst{
 &\dot{p}= -\gamma \MASST^{-1}p + \sigma \dot{W} - \Dq \xi \, \dot{\lambda}\COMMA &\\
 &p^{T} \MASST^{-1} \Dq \xi = 0\COMMA & \TAG{C_p}
}\end{equation}
with the usual fluctuation/dissipation relation $\sigma \sigma^T = 2\beta^{-1}\gamma$.
The Gaussian distribution of impulses
\be\label{e:maxwell}
\frac{1}{Z} \EXP{-\frac{\beta}{2} p^{T}\MASST^{-1}  p} \sigma_{T^{*}_{q}\mathcal{M}}(dp)
\ee
is invariant under the dynamics \VIZ{e:sdeimp}.
\bpro[Exact sampling of stochastic perturbation]\label{p:exactscheme}
Given the mass matrix $M$, suppose either $\dt$ or $\gamma$  are small enough so that the condition
\be
\frac{\delta t}{2}  \MASST^{-1}   \leq \gamma
\ee
holds in the sense of symmetric semi-definite matrices. Let $U$ be a centered and normalized Gaussian vector.
Consider the mid-point Euler scheme with constraints
\begin{equation} \label{e:midpoint}\syst{
& p_{n+1}= p_{n}  - \frac{\dt}{2} \gamma \MASST^{-1}(p_{n} + p_{n+1)} + \sqrt{\dt} \sigma U - \Dq \xi \,\lambda_{n+1} &\\
& p_{n+1}^T\MASST^{-1}  \Dq \xi = 0 \COMMA & \TAG{C_p} %
}\end{equation}
where $\lambda_{n+1}$ is the Lagrange multiplier associated with the constraint $(C_p)$.
The Markov kernel defined by the transition $p_{n}\to p_{n+1}$ is reversible with respect to
the Gaussian distribution \eqref{e:maxwell}.
\epro
\begin{proof}
After calculating the Lagrange multiplier the expression \VIZ{e:midpoint} can be written as
\[
p_{n+1} = p_{n}  - \frac{\dt}{2} P \gamma P^{T} \MASST^{-1}(p_{n} + p_{n+1}) + \sqrt{\dt} P \sigma U \PERIOD
\]
Consider the new variable
$\tilde{p} = \beta^{1/2} \MASST^{-1/2} p$,
and define the symmetric matrix
$$
\ZETAM \equiv \frac{\dt}{2} \MASST^{-1/2} P \gamma P^{T} \MASST^{-1/2} \COMMA
$$
as well as $\KAPPAM$, such that $\KAPPAM \KAPPAM^{T} = \ZETAM$.
In terms of these new variables we obtain from (\ref{e:midpoint})
\be \label{e:midpoint_tilde}
\tilde{p}_{n+1} = (\ID + \ZETAM)^{-1}(\ID-\ZETAM)\,\tilde{p}_n + 2 (\ID+\ZETAM)^{-1}\KAPPAM \,U \PERIOD
\ee
Moreover, the product measure $\sigma_{T^*_q\SMAN}(dp_n)\,\sigma_{T^*_q\SMAN}(dp_{n+1})$ is the measure induced on the linear subspace
of constraints by the scalar product $\MASST^{-1}$ and the Lebesgue measure $dp_n \, dp_{n+1}$.
Thus in the variables $( \tilde{p}_n, \tilde{p}_{n+1})$ this measure becomes, up to a constant, the measure induced by the usual Euclidean structure.
As a consequence the $\log$ density of the random variable $(\tilde{p}_n, \tilde{p}_{n+1})$ defined by
\eqref{e:midpoint_tilde} with respect to this latter measure is equal to
\begin{eqnarray*}
 && - \frac{1}{2} \abs{\tilde{p}_n}^2 - \frac{1}{8}\left(\tilde{p}_{n+1}- (\ID+\ZETAM)^{-1}(\ID-\ZETAM)\,\tilde{p}_n\right)^T
                  \ZETAM^{-1}(\ID+\ZETAM)^2 \left(\tilde{p}_{n+1}- (\ID+\ZETAM)^{-1}(\ID-\ZETAM)\,\tilde{p}_n \right) \\
&=& - \frac{1}{8}\tilde{p}_{n+1}^T \ZETAM^{-1}(\ID+\ZETAM)^2\,\tilde{p}_{n+1}
    - \frac{1}{8}\tilde{p}_{n}^T \ZETAM^{-1}(\ID+\ZETAM)^2\,\tilde{p}_{n} \COMMA
\end{eqnarray*}
which is indeed symmetric between $\tilde{p}_{n}$ and $\tilde{p}_{n+1}$. Hence
we have shown the reversibility of the induced Markov kernel and consequently
stationarity of the canonical Gaussian distribution.
\end{proof}

\bibliographystyle{plain}
\bibliography{mathias}

\begin{thebibliography}{10}

\bibitem{Arn89}
V.~I. Arnol'd.
\newblock {\em Mathematical methods of classical mechanics}.
\newblock Springer-Verlag, New York, 1989.

\bibitem{Beccaria94}
M.~Beccaria and G.~Curci.
\newblock The {K}ramers equation simulation algorithm: 1. {O}perator analysis.
\newblock {\em Phys. Rev. D}, 49:2578--2589, 1994.

\bibitem{Ben75}
C.H. Bennett.
\newblock Mass tensor molecular dynamics.
\newblock {\em J. Comp. Phys.}, 19:267--279, 1975.

\bibitem{BesStu}
A.~Beskos, G.O. Roberts, A.M. Stuart, and J.~Voss.
\newblock An {MCMC} method for diffusion bridges.
\newblock {\em Pre-print}, 2007.

\bibitem{BorSch97}
F.~Bornemann and C.~Schuette.
\newblock Homogenization of {H}amiltonian system with a strong constraining
  potential.
\newblock {\em Physica D}, 102:57--77, 1992.

\bibitem{LebLeg07}
C.~Le Bris and F.~Legoll.
\newblock Derivation of symplectic numerical schemes for highly oscillatory
  {H}amiltonian systems.
\newblock {\em C. R. Acad. Sci. Paris}, 344:277--282, 2007.

\bibitem{CanLeg05}
E.~Cances, F.~Legoll, and G.~Stoltz.
\newblock Comparison of {NVT} sampling methods.
\newblock Technical Report 2040, IMA, 2005.

\bibitem{CicLelVan06}
G.~Ciccotti, T.~Leli{\`e}vre, and E.~Vanden-Eijnden.
\newblock Projection of diffusions on submanifolds: Application to mean force
  computation.
\newblock Technical Report 309, CERMICS, 2006.

\bibitem{ClaBakStiBra94}
M.~E. Clamp, P.~G. Baker, C.J.Stirling, and A.Brass.
\newblock Hybrid {M}onte {C}arlo : an efficient algorithm for condensed matter
  simulation.
\newblock {\em J. Comput. Chem.}, 15(8):838--846, 1994.

\bibitem{JahLub06}
D.~Cohen, T.~Jahnke, K.~Lorenz, and Ch. Lubich.
\newblock Numerical integrators for highly oscillatory {H}amiltonian systems: a
  review.
\newblock In Aexander Mielke, editor, {\em Analysis, Modeling and Simulation of
  Multiscale Problems}, pages 553--576. Springer, 2006.

\bibitem{Creutz89}
M.~Creutz and A.~Gocksch.
\newblock Higher-order hybrid {M}onte {C}arlo algorithms.
\newblock {\em Phys. Rev. Lett.}, 63(1):9--12, 1988.

\bibitem{Duane85}
S.~Duane.
\newblock {Stochastic quantization versus the microcanonical ensemble - getting
  the best of both worlds}.
\newblock {\em {Nuclear Physics B}}, {257}({5}):{652--662}, {1985}.

\bibitem{DuaneKennedy87}
S.~Duane, A.~D. Kennedy, B.~J. Pendleton, and D.~Roweth.
\newblock Hybrid {M}onte {C}arlo.
\newblock {\em Phys. Lett.}, 195B(2):216--222, 1987.

\bibitem{HaiLub02}
C.~Lubich E.~Hairer, G.~Wanner.
\newblock {\em Geometric Numerical Integration, Structure-Preserving Algorithms
  for Ordinary Differential Equations}.
\newblock Springer, 2002.

\bibitem{CicVan06}
G.~Ciccotti E.~Vanden-Eijnden.
\newblock Second-order integrators for {L}angevin equations with holonomic
  constraints.
\newblock {\em Chem. Phys. Lett.}, 429(1-3):310--316, 2006.

\bibitem{EthKur86}
S.~N. Ethier and T.~G. Kurtz.
\newblock {\em Markov Processes: Characterization and Convergence}.
\newblock Wiley Series in Probability and Statistics, 1986.

\bibitem{Fed69}
H.~Federer.
\newblock {\em Geometric measure theory}.
\newblock Springer, New-Yoerk, 1969.

\bibitem{Fel71}
W.~Feller.
\newblock {\em An Introduction to Probability Theory and Its Applications}.
\newblock Wiley, New-York, 1971.

\bibitem{Fix78}
M.~Fixman.
\newblock Simulation of polymer dynamics. {I.} {G}eneral theory.
\newblock {\em J. Chem. Phys.}, 69:1527, 1978.

\bibitem{Fix86}
M.~Fixman.
\newblock Implicit algorithm for {B}rownian dynamics of polymers.
\newblock {\em Macromolecules}, 19:1195--1204, 1986.

\bibitem{CicKapVan05}
E.~Vanden-Eijnden G.~Ciccotti, R.~Kapral.
\newblock Blue moon sampling, vectorial reaction coordinates, and unbiased
  constrained dynamics.
\newblock {\em J. Chem. Phys.}, 6(9):1809--14, 2005.

\bibitem{Horowitz91}
A.~M. Horowitz.
\newblock A generalized guided {M}onte {C}arlo algorithm.
\newblock {\em Phys. Lett.B}, 268:247--252, 1991.

\bibitem{Kennedy98}
I.~Horv{\'{a}}th and A.~D. Kennedy.
\newblock The {L}ocal {H}ybrid {M}onte {C}arlo algorithm for free field theory.
\newblock {\em Nucl. Phys. B}, 510:367--400, 1998.

\bibitem{IzaHam04}
J.~A. Izaguirre and S.~S. Hampton.
\newblock Shadow hybrid {M}onte {C}arlo: an efficient propagator in phase space
  of macromolecules.
\newblock {\em J. Comput. Phys.}, 200(2):581--604, 2004.

\bibitem{RycCic77}
H.J.C.~Berendsen J-P. Ryckaert; G.~Ciccotti.
\newblock Numerical integration of the {C}artesian equations of motion of a
  system with constraints: Molecular dynamics of n-{A}lkanes.
\newblock {\em J. of Comp. Phys.}, 23:327--341, 1977.

\bibitem{FeeHesBer90}
H.J.C.~Berendsen K.~A.~Feenstra, B.~Hess.
\newblock Improving efficiency of large time-scale molecular dynamics
  simulations of hydrogen-rich systems.
\newblock {\em J. of Comp. Chemistry}, 20:786--798, 1999.

\bibitem{Kam85}
N.G.~Van Kampen.
\newblock Elimination of fast variables.
\newblock {\em Physics Reports}, 124(2):9--160, 1985.

\bibitem{LeiSke94}
B.~J. Leimkuhler and R.~D. Skeel.
\newblock Symplectic numerical integrators in constrained {H}amiltonian
  systems.
\newblock {\em J. Comput. Phys.}, 112(1):117--125, 1994.

\bibitem{MaoFrie90}
B.~Mao and A.R. Friedman.
\newblock Molecular dynamics simulation by atomic mass weighting.
\newblock {\em Biophysical Journal}, 58:803--805, 1990.

\bibitem{Met53}
N.~Metropolis, A.~W. Rosenbluth, M.~N. Rosenbluth, A.~H. Teller, and E.~Teller.
\newblock Equations of state calculations by fast computing machine.
\newblock {\em J. Chem. Phys.}, 21:1087--1091, 1953.

\bibitem{Osk92}
B.~Oksendal.
\newblock {\em Stochastic differential equations (3rd ed.): an introduction
  with applications}.
\newblock Springer-Verlag, 1992.

\bibitem{Rei95}
S.~Reich.
\newblock Smoothed dynamics of highly oscillatory {H}amiltonian systems.
\newblock {\em Physica D}, 89:28--42, 1995.

\bibitem{Rei00}
S.~Reich.
\newblock Smoothed {L}angevin dynamics of highly oscillatory systems.
\newblock {\em Physica D}, 138:210--224, 2000.

\bibitem{LeiRei05}
S.~Reich and B.~Leimkuhler.
\newblock {\em Simulating {H}amiltonian Dynamics}, volume~14 of {\em Cambridge
  Monographs on Applied and Computational Mathematics}.
\newblock Cambridge University Press, 2005.

\bibitem{RubUng57}
H.~Rubin and P.~Ungar.
\newblock Motion under a strong constraining force.
\newblock {\em Comm. Pure Appl. Math}, 10:65--87, 1957.

\bibitem{SchFis99}
Ch. Sch\"{u}tte, A.~Fischer, W.~Huisinga, and P.~Deuflhard.
\newblock A direct approach to conformational dynamics based on hybrid {M}onte
  {C}arlo.
\newblock {\em J. Comput. Phys.}, 151(1):146--168, 1999.

\bibitem{Tak80}
F.~Takens.
\newblock Motion under the influence of a strong constraining force.
\newblock In {\em Global theory of dynamical systems (Proc. Internat. Conf.,
  Northwestern Univ., Evanston, Ill., 1979)}, volume 819 of {\em Lecture Notes
  in Math.}, pages 425--445. Springer, Berlin, 1980.

\bibitem{GunBer77}
W.F. van Gunsteren and H.~J.~C. Berendsen.
\newblock Algorithms for macromolecular dynamics and constraint dynamics.
\newblock {\em Molecular Physics}, 34(5):1311--1327, 1977.

\bibitem{GunBer81}
W.F. van Gunsteren and H.J.C. Berendsen.
\newblock Algorithms for {B}rownian dynamics.
\newblock {\em Molecular Physics}, 45:637--647, 1982.

\end{thebibliography}

\end{document}